\documentclass[12pt]{article}

\usepackage{amsfonts,graphicx,color,epsfig,url}
\usepackage[letterpaper,asymmetric,left=1.1in,right=1.1in,top=1.25in,bottom=1.25in,
bindingoffset=0.0in]{geometry}

\title{Computing arithmetic invariants for hyperbolic reflection groups}

\author{Omar Antol\'{\i}n-Camarena\footnote{oantolin@math.utoronto.ca},
        Gregory R. Maloney\footnote{maloneyg@math.utoronto.ca},
        and Roland K. W. Roeder \footnote{rroeder@math.utoronto.ca}}
                                                                                    
\input{psfig.sty}
                                                                                    
\newtheorem{thm}{Theorem}

\newtheorem{lem}[thm]{Lemma}
\newtheorem{defn}[thm]{Definition}

\newcommand{\SO}{\ensuremath{\mathop{SO}\nolimits}}
\newcommand{\SL}{\ensuremath{\mathop{SL}\nolimits}}

\newcommand{\PSL}{\ensuremath{\mathop{PSL}\nolimits}}
\newcommand{\Isom}{\ensuremath{\mathop{\mathrm{Isom}}\nolimits}}
\newcommand{\id}{\ensuremath{\mathop{\mathrm{id}}\nolimits}}
\newcommand{\tr}{\ensuremath{\mathop{\mathrm{tr}}\nolimits}}
\newcommand{\hilb}[3]{\ensuremath{\left( \frac{#1, #2}{#3} \right)}}

\newcommand{\ring}[1]{\ensuremath{R_{#1}}}
                                                                                    
\begin{document}

\maketitle

\abstract{ 
We describe a collection of computer scripts written in PARI/GP to compute, for
reflection groups determined by finite-volume polyhedra in $\mathbb{H}^3$, the
commensurability invariants known as the invariant trace field and invariant
quaternion algebra.  Our scripts also allow one to determine arithmeticity of such
groups and the isomorphism class of the invariant quaternion algebra by
analyzing its ramification.  

We present many computed examples of these invariants.  This is enough to show
that most of the groups that we consider are pairwise incommensurable.  For
pairs of groups with identical invariants, not all is lost: when both groups
are arithmetic, having identical invariants guarantees commensurability.  We
discover many ``unexpected'' commensurable pairs this way.  We also present a
non-arithmetic pair with identical invariants for which we cannot determine
commensurability.

}

\vspace{.1in}

\section{Introduction}

Suppose that $P$ is a finite-volume polyhedron in $\mathbb{H}^3$ each of whose
dihedral angles is an integer submultiple of $\pi$.  Then the group
$\Lambda(P)$ generated by reflections in the faces of $P$ is a discrete subgroup
of $\Isom(\mathbb{H}^3)$.  If one restricts attention to the subgroup
$\Gamma(P)$ consisting of orientation-preserving elements of $\Lambda(P)$, one
naturally obtains a discrete subgroup of $\PSL(2,\mathbb{C}) \cong
\Isom^+(\mathbb{H}^3)$.  This very classical family of
finite-covolume Kleinian groups is known as the family of {\em polyhedral reflection
groups}.  

There is a complete classification of hyperbolic polyhedra with non-obtuse
dihedral angles, and hence of hyperbolic reflection groups, given by Andreev's
Theorem \cite{AND,ROE2} (see also \cite{RH,H,STEVE} for alternatives to the
classical proof); however, many more detailed questions about the resulting
reflection group remain mysterious.   We will refer to finite-volume hyperbolic
polyhedra with non-obtuse dihedral angles as {\em Andreev Polyhedra} and
finite-volume hyperbolic polyhedra with integer submultiple of $\pi$ dihedral
angles {\em Coxeter Polyhedra}.

A fundamental question for general Kleinian groups is: given $\Gamma_1$ and
$\Gamma_2$ does there exist an appropriate conjugating element $g\in
\PSL(2,\mathbb{C})$ so that $\Gamma_1$ and $g\Gamma_2g^{-1}$ both have a
finite-index subgroup in common?  In this case, $\Gamma_1$ and $\Gamma_2$ are
called {\em commensurable}.  Commensurable Kleinian groups have many properties
in common, including coincidences in the lengths of closed geodesics (in the
corresponding orbifolds) and a rational relationship (a commensurability)
between their covolumes, if the groups are of finite-covolume.  See
\cite{NEUMANN_REID} for many more interesting aspects of commensurability in
the context of Kleinian groups.

\begin{figure}
\begin{center}
\includegraphics[scale=0.9]{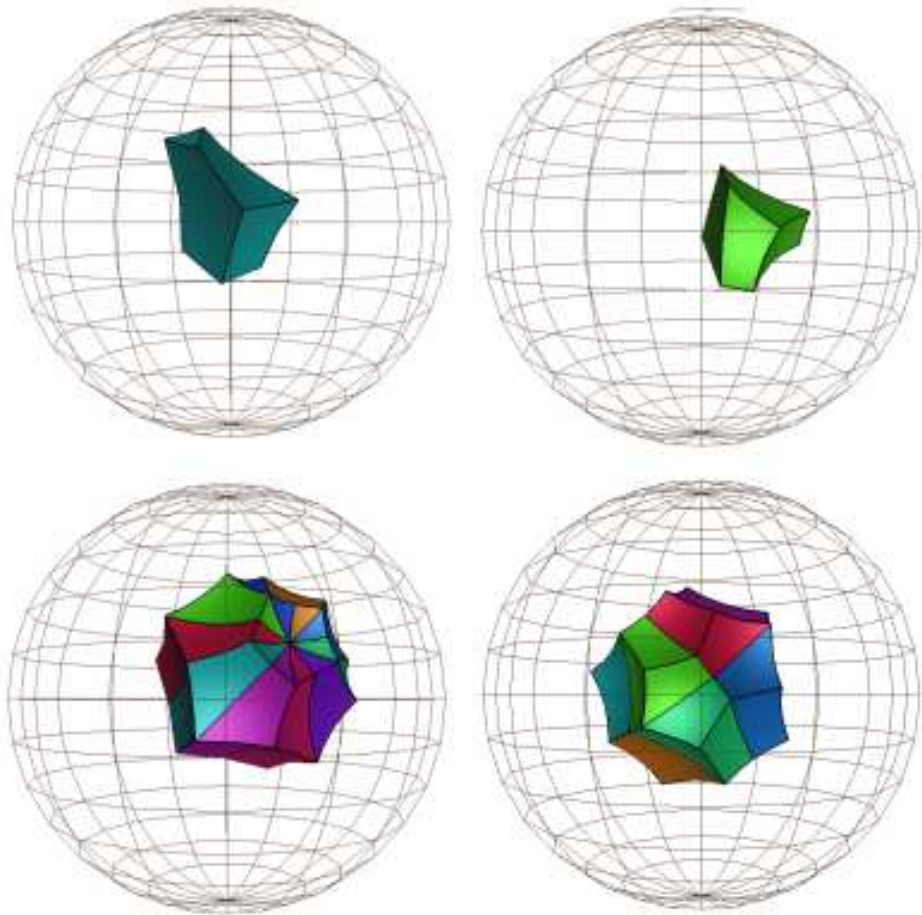}
\end{center}
\caption{\label{FIG:COMMENS} Two commensurable polyhedra $P_1$ (left top) and $P_2$ (right top)
which tile a common larger polyhedron $Q$, here the right-angled dodecahedron.}  
\end{figure}

If $\Gamma_1$ and $\Gamma_2$ are fundamental groups of hyperbolic manifolds
$M_1$ and $M_2$, commensurability is the same as the existence of a common
finite-sheeted cover $\widetilde{M}$ of $M_1$ and of $M_2$.  Similarly, if
$\Gamma(P_1)$ and $\Gamma(P_2)$ are polyhedral reflection groups, they are
commensurable if and only if there is a larger polyhedron $Q$ that is tiled
both by $P_1$ under reflections in the faces (of $P_1$) and by $P_2$ under
reflections in the faces (of $P_2$.)   The existence of such a polyhedron $Q$
is clearly a fundamental and delicate question from hyperbolic geometry.  See
Figure \ref{FIG:COMMENS} for an example of two commensurable polyhedra.  These
coordinates for these polyhedra were computed using \cite{ROE_CONSTRUCTION} and
displayed in the conformal ball model of $\mathbb{H}^3$ using Geomview
\cite{GEOMVIEW}.

A pair of sophisticated number-theoretic invariants has been developed by
Reid and others to distinguish between commensurability
classes of general finite-covolume Kleinian groups.  See the recent textbook
\cite{MR_BOOK} and the many references therein.  Given $\Gamma$, these
invariants are a number field $k(\Gamma)$ known as the {\em invariant trace
field} and a quaternion algebra $A(\Gamma)$ over $k(\Gamma)$ known as the {\em
invariant quaternion algebra}.  In fact, the invariant trace field is obtained
by intersecting all of the fields generated by traces of elements of the
finite-index subgroups of $\Gamma$.  It is no surprise that such a field is
related to commensurability because the trace of a loxodromic element $a$ of
$\Gamma$ is related to the translation distance $d$ along the axis of $a$ by
$2\cosh(d) = {\rm Re}(tr(A))$.

The pair $\left(k(\Gamma),A(\Gamma)\right)$ does a pretty good job to
distinguish commensurability classes, but there are examples of incommensurable
Kleinian groups with the same $\left(k(\Gamma),A(\Gamma)\right)$.  For
arithmetic groups, however, the pair $\left(k(\Gamma),A(\Gamma)\right)$ is a
complete commensurability invariant.  Thus, one can find unexpected
commensurable pairs of groups by computing these two invariants and by
verifying that each group is arithmetic.  See Subsection
\ref{SUBSEC:COMMENS_PAIRS} for examples of such pairs that were discovered in
this way.  The precise definitions of the invariant trace field, the invariant
quaternion algebra, and an arithmetic group will be given in Section
\ref{SEC:ITF}.

It can be rather difficult to compute the invariant trace field and invariant
quaternion algebra of a given Kleinian group ``by hand.''  However there is a
beautiful computer program called SNAP \cite{SNAP} written by Coulson, Goodman,
Hodgson, and Neumann, as described in \cite{SNAP_PAPER}.  They have computed
the invariant trace field and invariant quaternion algebra, as well as many other
interesting invariants, for many of the manifolds in the Hildebrand-Weeks
census \cite{HILDEBRAND_WEEKS} and in the Hodgson-Weeks census
\cite{HODG_WEEKS}.  The basic idea used in SNAP is to compute a high-precision
decimal approximation for an ideal triangulation of the desired manifold $M$
using Newton's Method and then to use the LLL algorithm \cite{LLL} to guess
exact algebraic numbers from the approximate values.  These guessed values can
be checked for correctness using the gluing equations describing $M$, and if
the values are correct, the invariant trace field and invariant quaternion algebra
can be computed from the exact triangulation. 

SNAP provides a vast source of examples, also seen
in the appendix of the book \cite{MR_BOOK}, and adds enormous flavor to the field.
The fundamental techniques used in SNAP provide inspiration for our
current work with polyhedral reflection groups.

In the case of polyhedral reflection groups there is a simplified description
of the invariant trace field and the invariant quaternion algebra in terms of
the Gram matrix of the polyhedron \cite{MR_POLYHEDRA}.  This theorem avoids the
rather tedious trace calculations and manipulation of explicit generators of the
group.  Following the general technique used in the program SNAP we compute  a
set of outward unit normals to the faces of the polyhedron $P$ to a high
decimal precision and then use the LLL algorithm to guess the exact normals as
algebraic numbers.  From these normals the Gram matrix is readily computed,
both allowing us to check whether the guessed algebraic numbers are in fact
correct, and providing the exact data needed to use the theorem from
\cite{MR_POLYHEDRA} in order to compute the invariant trace field and invariant
quaternion algebra for $\Gamma(P)$. 

Our technique is illustrated for a simple example in Section
\ref{SEC:WORKED_EXAMPLE} and a description of our program (available to
download, see \cite{SNAPHEDRON}) is given in Section \ref{SEC:PROGRAM}.  Section
\ref{SEC:QUATERNION} provides details on how to interpret the quaternion
algebra.  Finally in Section \ref{SEC:RESULTS} we provide many results of our
computations.

\vspace{0.15in}
\noindent
{\center \bf Acknowledgments:}
We thank Colin Maclachlan and Alan Reid for their beautiful work and exposition
on the subject in \cite{MR_BOOK}.  We also thank Alan Reid for his many helpful
comments.

We thank the authors of SNAP \cite{SNAP} and the corresponding paper
\cite{SNAP_PAPER}, which inspired this project (including the choice of name
for our collection of scripts).  Among them Craig Hodgson has provided helpful
comments.  

We thank Andrei Vesnin who informed us about his result about arithmeticity of
L\"obell polyhedra.

We effusively thank the writers of PARI/GP \cite{PARI2}, the system
in which we have written our entire program and which is also used in SNAP
\cite{SNAP}.

The third author thanks Mikhail Lyubich and Ilia Binder for their financial
support and interest in the project.  He also thanks John Hubbard, to whom this
volume is dedicated, for introducing him to hyperbolic geometry and for his
enthusiasm for mathematics in general and experimental mathematics in
particular.

\section{Hyperbolic polyhedra and the Gram Matrix}

We briefly recall some fundamental hyperbolic geometry, including the
definition of a hyperbolic polyhedron and of the Gram matrix of a polyhedron.

Let $E^{3,1}$ be the four-dimensional Euclidean space with the indefinite
metric $\Vert {\bf x} \Vert^2 = -x_0^2+x_1^2+x_2^2+x_3^2$.  Then hyperbolic
space $\mathbb{H}^3$ is the component having $x_0 > 0$ of the subset of $E^{3,1}$ given by

$$\Vert {\bf x} \Vert^2 = -x_0^2+x_1^2+x_2^2+x_3^2 = -1$$

\noindent
with the Riemannian metric induced by the indefinite
metric

$$-dx_0^2+dx_1^2+dx_2^2+dx_3^2.$$

\vspace{.05in}

The hyperplane orthogonal to a vector ${\bf v} \in E^{3,1}$ intersects
$\mathbb{H}^3$ if and only if $\langle{\bf v},{\bf v}\rangle> 0$.  Let ${\bf v}
\in E^{3,1}$ be a vector with $\langle{\bf v},{\bf v}\rangle > 0$, and define

$$P_{\bf v} = \{{\bf w} \in \mathbb{H}^3 | \langle{\bf w},{\bf v}\rangle
= 0\}$$

\noindent to be the hyperbolic plane orthogonal to ${\bf v}$;  and the
corresponding closed half space:

$$H_{\bf v}^+ = \{{\bf w} \in \mathbb{H}^3 | \langle{\bf w},{\bf v}\rangle 
\geq 0 \}.$$

\noindent 
It is a well known fact that
given two planes $P_{\bf v}$ and $P_{\bf w}$ in 
$\mathbb{H}^3$ with $\langle{\bf v},{\bf v}\rangle = 1$ and $\langle{\bf 
w},{\bf w}\rangle = 1$, they:

\begin{itemize}
\item intersect in a line if and only if $\langle{\bf v},{\bf w}\rangle^2
< 1$, in which case their dihedral angle is $\arccos(-\langle{\bf v},{\bf
w}\rangle)$.

\item intersect in a single point at infinity if and only if $\langle{\bf
v},{\bf w}\rangle^2 = 1$, in this case their dihedral angle is $0$.

\item are disjoint if and only if $\langle{\bf v},{\bf w}\rangle^2 > 1$, in
which  case the distance between them is ${\rm arccosh}(-\langle{\bf v},{\bf
w}\rangle))$.

\end{itemize}

Suppose that ${\bf e}_1,\ldots,{\bf e}_n$ satisfy $\langle {\bf e}_i,{\bf e}_i
\rangle > 0$ for each $i$.  Then, a {\it hyperbolic polyhedron} is an
intersection

$$P = \bigcap_{i=0}^n H_{{\bf e}_i}^+ $$

\noindent having non-empty interior.  

If we normalize the vectors ${\bf e}_i$ that are orthogonal to the faces of a
polyhedron $P$, the {\em Gram Matrix} of $P$ is given by $M_{ij}(P) =2\langle
{\bf e}_i, {\bf e}_j \rangle$.  It is also common to define the Gram matrix
without this factor of $2$, but our definition is more convenient for arithmetic
reasons.  By construction, a Gram matrix is symmetric and has $2$s on the
diagonal.  Notice that the Gram matrix encodes information about both the
dihedral angles between adjacent faces of $P$ and the hyperbolic
distances between non-adjacent faces.

\section{Invariant Trace Field, Invariant Quaternion Algebra, and
Arithmeticity}\label{SEC:ITF}

The \emph{trace field} of a subgroup $\Gamma$ of $\PSL (2,
\mathbb{C})$ is the field generated by the traces of its
elements; that is, $\mathbb{Q}(\tr \Gamma) := \mathbb{Q}(\tr
\gamma : \gamma \in \Gamma)$.{\footnote{Note that for $\gamma \in \PSL
    (2, \mathbb{C})$, the trace $\tr \gamma$ is only defined up to
    sign.}} This field is not a commensurablitity invariant as shown
by the following example found in {\cite{MR_BOOK}}.

Consider the group $\Gamma$ generated by
\[ A = \left(\begin{array}{cc}
     1 & 1\\
     1 & 0
   \end{array}\right), \hspace{1em} B = \left(\begin{array}{cc}
     1 & 0\\
     - \omega & 1
   \end{array}\right), \]
where $\omega = (- 1 + i \sqrt{3}) / 2$. The trace field of $\Gamma$ is
$\mathbb{Q}( \sqrt{- 3})$. Now let $X = \left(\begin{array}{cc}
  i & 0\\
  0 & - i
\end{array}\right)$. It is easy to see that $X$ normalizes $\Gamma$ and its
square is the identity (in $\PSL (2, \mathbb{C})$), so that $\Lambda =
\left\langle \Gamma, X \right\rangle$ contains $\Gamma$ as a subgroup
of index $2$ and is therefore commensurable with $\Gamma$. But
$\Lambda$ also contains $XBA = \left(\begin{array}{cc}
    i & i\\
    i \omega &  -i + i\omega
\end{array}\right)$, so the trace field of $\Lambda$ contains $i$ in addition
to $\omega$.

The easiest way to fix this, that is, to get a commensurability
invariant related to the trace field, is to associate to $\Gamma$ the
intersection of the trace fields of all finite index subgroups of
$\Gamma$; this is the \emph{invariant trace field} denoted $k \Gamma$.

While this definition clearly shows commensurability invariance, it
does not lend itself to practical calculation. The proof of Theorem
3.3.4 in {\cite{MR_BOOK}} brings us closer: it shows that instead of
intersecting many trace fields, one can look at a single one, namely,
the invariant trace field of $\Gamma$ equals the trace field of its
subgroup $\Gamma^{(2)} : = \langle \gamma^2 : \gamma \in \Gamma
\rangle$. (This also shows that the invariant trace field is non-trivial which is not entirely clear from the definition as an intersection.)  When a finite set of generators for the group is known, this
is actually enough to compute the invariant trace field. Indeed,
Lemma 3.5.3 in \cite{MR_BOOK} establishes that if $\Gamma = \langle
\gamma_1, \gamma_2, \ldots, \gamma_n \rangle$, the invariant trace
field of $\Gamma$ is generated by $\{\tr(\gamma_i) : 1\le i \le n\}
\cup \{\tr(\gamma_i \gamma_j) : 1\le i < j \le n\} \cup \{\tr(\gamma_i
\gamma_j \gamma_k : 1\le i < j < k \le n\}$.

For reflection groups a more efficient description can be given in terms of the
Gram matrix. The description given above uses roughly $n^3/6$ generators for a
polyhedron with $n$ faces; the following description will use only around
$n^2/2$. But before we state it we need to define a certain quadratic space
over the field $k (P):=\mathbb{Q}(a_{i_1 i_2} a_{i_2 i_3} \cdots a_{i_r i_1} :
\{i_1, i_2, \ldots, i_r \} \subset \{1, 2, \ldots, n\})$ associated to a
polyhedron; this space will also appear in the next section in the theorem used
to calculate the invariant quaternion algebra.

As in the previous section, given a polyhedron $P$ we will denote the
outward-pointing normals to the faces by ${\bf e}_1,\ldots,{\bf e}_n$
and the Gram matrix by $(a_{ij})$. Define $M(P)$ as the vector space over $k(P)$ spanned by of all the vectors of
the form $a_{1 i_1} a_{i_1 i_2} \cdots a_{i_{r - 1} i_r} {\bf
  e}_{i_r}$ where $\{i_1, i_2, \ldots, i_r\}$ ranges over the subsets
of $\{1,2,\ldots,n\}$ and $n$ is the number of faces of $P$.  This
space $M(P)$ will be equipped with the restriction of the quadratic
form with signature $(3,1)$ used in $\mathbb{H}^3$. We recall that
the discriminant of a non-degenerate symmetric bilinear form $\langle
\cdot, \cdot \rangle$ is defined as $\det \left(\langle v_i, v_j
  \rangle \right)_{ij}$ where $\{v_i\}_i$ is a basis for the vector
space on which the form is defined. The discriminant does depend on
the choice of basis, but for different bases the discriminants differ
by multiplication by a square in the ground field: indeed, if $u_i =
\sum_j \alpha_{ij} v_j$, the discriminant for the basis $\{u_i\}_i$ is
that of the basis $\{v_i\}_i$ multiplied by $\det (\alpha_{ij})^2$.

Now we can state the theorem we use to calculate the invariant trace
field, Theorem 10.4.1 in {\cite{MR_BOOK}}:

\begin{thm}
  \label{THM:POLY_ITF}
Let $P$ be a Coxeter polyhedron 
and let $\Gamma$ be the reflection group it determines.
Let $(a_{ij})$ be the Gram matrix of $P$. The invariant trace field of $\Gamma$
is $k (P) ( \sqrt{d})$, where $d$ is the discriminant of the quadratic space
$M(P)$ and $k (P) =\mathbb{Q}(a_{i_1 i_2} a_{i_2 i_3} \cdots a_{i_r i_1} :
\{i_1, i_2, \ldots, i_r \} \subset \{1, 2, \ldots, n\})$ is the field defined
previously.
\end{thm}

\subsection{The Invariant Quaternion Algebra}

A \emph{quaternion algebra} over a field $F$ is a four-dimensional
associative algebra $A$ with basis $\{1, i, j, k\}$ satisfying $i^2 =
a 1$, $j^2 = b 1$ and $ij = ji = - k$ for some $a, b \in F$. Note that
$k^2 = (ij)^2 = ijij = - ijji = - ab 1$. The case $F =\mathbb{R}$, $a
= b = - 1$ gives Hamilton's quaternions.

The quaternion algebra defined by a pair $a, b$ of elements of $F$ is denoted
by its \emph{Hilbert symbol} $\hilb{a}{b}{F}$.  A quaternion algebra does not
uniquely determine a Hilbert symbol for it, since, for example, $\hilb{a}{b}{F}
= \hilb{a}{-ab}{F} = \hilb{au^2}{bv^2}{F}$ for any invertible elements $u, v
\in F$. Fortunately, there is a computationally effective way of deciding
whether two Hilbert symbols give the same quaternion algebra. This will be
discussed in Section \ref{SEC:QUATERNION}; for now we will just define the
invariant quaternion algebra of a subgroup of $\PSL (2, \mathbb{C})$ and state
the theorem we use to calculate a Hilbert symbol for it.

Given any non-elementary{\footnote{This means that the action of
    $\Gamma$ on $\mathbb{H}^3 \cup \hat{\mathbb{C}}$ has no
    finite orbits. Reflection groups determined by finite-volume
    polyhedra are always non-elementary.}} subgroup $\Gamma$ of $\PSL
(2, \mathbb{C})$, we can form the algebra $A_0 \Gamma := \{ \sum a_i
\gamma_i : a_i \in \mathbb{Q}(\tr \Gamma), \gamma_i \in \Gamma \}$.
(Abusing notation slightly we consider the elements of $\Gamma$ as
matrices defined up to sign.)  This turns out to be a quaternion
algebra over the trace field $\mathbb{Q}(\tr \Gamma)$ (see Theorem
3.2.1 in {\cite{MR_BOOK}}).

Just as with the trace fields, we define the \emph{invariant
  quaternion algebra} of $\Gamma$, denoted by $A \Gamma$, as the
intersection of all the quaternion algebras associated to finite-index
subgroups of $\Gamma$.

When $\Gamma$ is finitely generated in addition to non-elementary, we
are in a situation similar to that of the invariant trace field in
that the invariant quaternion algebra is simply the quaternion algebra
associated to the subgroup $\Gamma^{(2)}$ of $\Gamma$, or in symbols, $A
\Gamma = A_0 \Gamma^{(2)}$. To see this, note that Theorem 3.3.5 in
{\cite{MR_BOOK}} states that for finitely generated non-elementary
$\Gamma$, the quaternion algebra $A_0 \Gamma^{(2)}$ is a
commensurability invariant. Now, given an arbitrary finite-index
subgroup $\Lambda$ of $\Gamma$ we have $A \Gamma \subset A_0
\Gamma^{(2)} = A_0 \Lambda^{(2)} \subset A_0 \Lambda$.

In the case where $\Gamma$ is the reflection group of a polyhedron
$P$, $A \Gamma$ can be identified as the even-degree subalgebra of a
certain Clifford algebra associated with $P$. Let us briefly recall
the basic notions related to Clifford algebras. Given an
$n$-dimensional vector space $V$ over a field $F$ equipped with a
non-degenerate symmetric bilinear form $\left\langle \cdot, \cdot
\right\rangle$ and associated quadratic form $\left\| \cdot
\right\|^2$, the Clifford algebra it determines is the
$2^n$-dimensional $F$-algebra generated by all formal products of
vectors in $V$ subject to the condition $v^2 = \left\langle v, v
\right\rangle 1$ (where $1$ is the empty product of vectors). If
$\{v_1, v_2, \ldots, v_n \}$ is an orthogonal basis of $V$, a basis
for the Clifford algebra $C (V)$ is $\{v_{i_1} v_{i_2} \cdots v_{i_r}
: 1 \leq r \leq n, 1 \leq i_1 < i_2 < \cdots < i_r \leq n\}$. There is
a $\mathbb{Z}_2$-grading of $C (V)$ given on monomials by the parity
of the number of vector factors.

Now we can state the result mentioned above: $A \Gamma$ is the even-degree
subalgebra of $C (M)$ where $M$ is the vector space over the invariant trace
field of $\Gamma$ that appears in Theorem \ref{THM:POLY_ITF}. This relationship
between the invariant trace algebra and $M$ allows one to prove a theorem
giving an algorithm for computing the Hilbert symbol for the invariant
quaternion algebra, part of Theorem 3.1 in {\cite{MR_POLYHEDRA}}:

\begin{thm}\label{THM:POLY_QA}
  The invariant quaternion algebra of the reflection group $\Gamma$ of
  a polyhedron $P$ is given by
  \[ A \Gamma = \hilb{ - \left\| u_1 \right\|^2 \left\| u_2
      \right\|^2}{- \left\| u_1 \right\|^2 \left\| u_3 \right\|^2}{k
      \Gamma}, \] where $\{u_1, u_2, u_3, u_4 \}$ is an
  orthogonal basis for the quadratic space $M(P)$ defined in the
  previous section.
\end{thm}

\vspace{0.1in}

In many cases when studying an orbifold $O = \mathbb{H}^3/\Gamma$, a simple
observation about the singular locus of $O$ leads to the fact that the invariant
quaternion algebra can represented by the Hilbert symbol $\hilb{ -1}{-1}{k
\Gamma}$.  This happens particularly often for polyhedral reflection groups.

\begin{figure}
\begin{center}
\begin{picture}(0,0)%
\epsfig{file=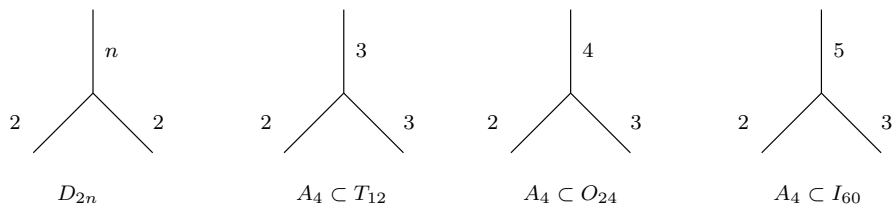}%
\end{picture}%
\setlength{\unitlength}{3947sp}%
\begingroup\makeatletter\ifx\SetFigFont\undefined%
\gdef\SetFigFont#1#2#3#4#5{%
  \reset@font\fontsize{#1}{#2pt}%
  \fontfamily{#3}\fontseries{#4}\fontshape{#5}%
  \selectfont}%
\fi\endgroup%
\begin{picture}(5659,1243)(676,-1667)
\put(2476,-1636){\makebox(0,0)[lb]{\smash{{\SetFigFont{8}{9.6}{\familydefault}{\mddefault}{\updefault}{\color[rgb]{0,0,0}$A_4 \subset T_{12}$}%
}}}}
\put(1576,-1186){\makebox(0,0)[lb]{\smash{{\SetFigFont{8}{9.6}{\familydefault}{\mddefault}{\updefault}{\color[rgb]{0,0,0}$2$}%
}}}}
\put(676,-1186){\makebox(0,0)[lb]{\smash{{\SetFigFont{8}{9.6}{\familydefault}{\mddefault}{\updefault}{\color[rgb]{0,0,0}$2$}%
}}}}
\put(1276,-736){\makebox(0,0)[lb]{\smash{{\SetFigFont{8}{9.6}{\familydefault}{\mddefault}{\updefault}{\color[rgb]{0,0,0}$n$}%
}}}}
\put(3151,-1186){\makebox(0,0)[lb]{\smash{{\SetFigFont{8}{9.6}{\familydefault}{\mddefault}{\updefault}{\color[rgb]{0,0,0}$3$}%
}}}}
\put(2251,-1186){\makebox(0,0)[lb]{\smash{{\SetFigFont{8}{9.6}{\familydefault}{\mddefault}{\updefault}{\color[rgb]{0,0,0}$2$}%
}}}}
\put(2851,-736){\makebox(0,0)[lb]{\smash{{\SetFigFont{8}{9.6}{\familydefault}{\mddefault}{\updefault}{\color[rgb]{0,0,0}$3$}%
}}}}
\put(4576,-1186){\makebox(0,0)[lb]{\smash{{\SetFigFont{8}{9.6}{\familydefault}{\mddefault}{\updefault}{\color[rgb]{0,0,0}$3$}%
}}}}
\put(3676,-1186){\makebox(0,0)[lb]{\smash{{\SetFigFont{8}{9.6}{\familydefault}{\mddefault}{\updefault}{\color[rgb]{0,0,0}$2$}%
}}}}
\put(4276,-736){\makebox(0,0)[lb]{\smash{{\SetFigFont{8}{9.6}{\familydefault}{\mddefault}{\updefault}{\color[rgb]{0,0,0}$4$}%
}}}}
\put(6151,-1186){\makebox(0,0)[lb]{\smash{{\SetFigFont{8}{9.6}{\familydefault}{\mddefault}{\updefault}{\color[rgb]{0,0,0}$3$}%
}}}}
\put(5251,-1186){\makebox(0,0)[lb]{\smash{{\SetFigFont{8}{9.6}{\familydefault}{\mddefault}{\updefault}{\color[rgb]{0,0,0}$2$}%
}}}}
\put(5851,-736){\makebox(0,0)[lb]{\smash{{\SetFigFont{8}{9.6}{\familydefault}{\mddefault}{\updefault}{\color[rgb]{0,0,0}$5$}%
}}}}
\put(3901,-1636){\makebox(0,0)[lb]{\smash{{\SetFigFont{8}{9.6}{\familydefault}{\mddefault}{\updefault}{\color[rgb]{0,0,0}$A_4 \subset O_{24}$}%
}}}}
\put(5476,-1636){\makebox(0,0)[lb]{\smash{{\SetFigFont{8}{9.6}{\familydefault}{\mddefault}{\updefault}{\color[rgb]{0,0,0}$A_4 \subset I_{60}$}%
}}}}
\put(976,-1636){\makebox(0,0)[lb]{\smash{{\SetFigFont{8}{9.6}{\familydefault}{\mddefault}{\updefault}{\color[rgb]{0,0,0}$D_{2n}$}%
}}}}
\end{picture}%
\end{center}
\caption{\label{FIG:SING} The possible vertices in the singular locus of a $3$-dimensional orbifold.  The last three have vertex stabilizer containing $A_4$.}
\end{figure}

Any vertex in the singular locus must be trivalent and must have labels from
the short list shown in Figure \ref{FIG:SING}.  See \cite{BP} for more information on orbifolds, and in particular page 24 from which Figure \ref{FIG:SING}
is essentially copied.

The last three have vertex stabilizer containing $A_4$, so if
$\mathbb{H}^3/\Gamma$ has singular locus containing such a vertex, $\Gamma$ must
contain $A_4$ as a subgroup.  In this case, the invariant quaternion algebra
can be represented by the Hilbert symbol $\hilb{ -1}{-1}{k \Gamma}$, see
\cite{MR_BOOK}, Lemma 5.4.1.  (See also Lemma 5.4.2.)

For a polyhedral reflection group generated by a Coxeter polyhedron $P$, the
corresponding orbifold $\mathbb{H}^3/\Gamma$ has underlying space
$\mathbb{S}^3$ and the singular set is (an unknotted) copy of the edge graph of
$P$.  The label at each edge having dihedral angle $\frac{\pi}{n}$ is merely
$n$.  In many cases, this simplification makes it easier to compute the Hilbert
symbol of a polyhedral reflection group.  We do not automate this check within
our program, but it can be useful to the reader.

\section{Quaternion algebras and their invariants}\label{SEC:QUATERNION} As
mentioned before, a quaternion algebra is fully determined by its Hilbert
symbol $\hilb{a}{b}{F}$, although this is by no means unique. For example,
\[\hilb{b}{a}{F},\qquad
\hilb{a}{-ab}{F},\qquad\textrm{and}\qquad\hilb{ax^2}{by^2}{F}\] all determine
the same algebra (here $x$ and $y$ are arbitrary invertible elements of $F$.)

Taking $F =\mathbb{R}$, the Hilbert symbol $\hilb{-1}{-1}{\mathbb{R}}$
represents the ordinary quaternions (or Hamiltonians,) denoted by
$\mathcal{H}$. If $F$ is any field, the Hilbert symbol $\hilb{1}{1}{F}$ is
isomorphic to $M_2(F)$, the two-by-two matrices over $F$.

A natural question now arises; namely, when do two Hilbert symbols represent
the same quaternion algebra? This question is pertinent for us---especially in
the case when $F$ is a number field---because the invariant quaternion algebra
is an invariant for a 
reflection group. For these purposes, we will need some way of classifying
quaternion algebras over number fields. All of the following material on
quaternion algebras appears in the reference \cite{MR_BOOK}.

The first step in the classification of quaternion algebras is the following
theorem: \begin{thm}\label{split}Let $A$ be a quaternion algebra over a field
$F$. Then either $A$ is a division algebra or $A$ is isomorphic to
$M_2(F)$.\end{thm} In the latter case, we say that $A$ \emph{splits}.
There are several different ways of expressing this condition, one of which
will be particularly useful for us: 

\begin{thm}\label{THM:HILBERT_EQ}
The
quaternion algebra $A = \hilb{a}{b}{F}$ splits over $F$ if and only if the
equation \begin{equation}\label{hilberteq}ax^2 + by^2 = 1\end{equation} has a
solution in $F\times F$. We call this equation the \emph{Hilbert equation} of
$A$.\end{thm}

When $F$ is a number field, it turns out that in order to classify the
quaternion algebras over $F$ completely we need to look at quaternion algebras
over the completions of $F$ with respect to its valuations. The first chapter
of \cite{MR_BOOK} contains a brief introduction to number fields and
valuations, and \cite{LANG_NUMBER} is a standard text on the subject.
\begin{defn}Let $F$ be any field. A \emph{valuation} on $F$ is a map $\nu :
F\rightarrow \mathbb{R}^+$ such that \begin{itemize} \item[(i)] $\nu (x) \geq
0$ for all $x\in F$ and $\nu (x) = 0$ if and only if $x = 0$, \item[(ii)] $\nu
(xy) = \nu (x)\nu(y)$ for all $x, y\in F$,\qquad\textrm{and} \item[(iii)] $\nu
(x + y)\leq \nu (x) + \nu (y)$ for all $x, y\in F$.  \end{itemize}\end{defn}

Any field admits a trivial valuation $\nu (x) = 1$ for all $x\neq 0$. When $F$
is a subfield of the real (or complex) numbers, the ordinary absolute value (or
modulus) function is a valuation when restricted to $F$. In general, valuations
on a field fall into two different classes.

\begin{defn}If a valuation $\nu$ on a field $F$ also satisfies
\begin{itemize}
\item[(iv)] $\nu (x + y) \leq max\{ \nu (x), \nu (y)\}$ for all $x, y\in F$,
\end{itemize}
then $\nu$ is called a \emph{non-Archimedean} valuation. If the valuation $\nu$ does not satisfy (iv), then it is called \emph{Archimedean}.
\end{defn}

There is also a notion of equivalence between valuations.

\begin{defn}Two valuations $\nu_1$ and $\nu_2$ on $F$ are called \emph{equivalent} if there exists some $\alpha\in\mathbb{R}^+$ such that $\nu_2(x) = \big( \nu_1(x)\big)^\alpha$ for all $x\in F$.
\end{defn}

When $F$ is a number field, it is possible to classify all valuations on $F$ up
to this notion of equivalence. Let $\sigma$ be a real or complex embedding of
$F$. Then a valuation $\nu_\sigma$ can be defined by $\nu_\sigma (x) = |\sigma
(x)|$, where $|\cdot |$ is the absolute value on $\mathbb{R}$ or modulus on
$\mathbb{C}$. This is an Archimedean valuation on $F$, and up to equivalence
these are the only Archimedean valuations that $F$ admits.

Denote by $\ring{F}$ the ring of integers of $F$---i.e. the set of elements of
$F$ satisfying some monic polynomial equation with integer coefficients---which
is a subring of $F$. Let $\mathcal{P}$ be a prime ideal in $\ring{F}$. Define a
function $n_\mathcal{P}:\ring{F}\rightarrow \mathbb{Z}$ by $n_\mathcal{P}(a) =
m$, where $m$ is the largest integer such that $a\in\mathcal{P}^m$. Since $F$
is the field of fractions of $\ring{F}$, $n_\mathcal{P}$ can be extended to all
of $F$ by the rule $n_\mathcal{P}(a/b) = n_\mathcal{P}(a)-n_\mathcal{P}(b)$.
Now pick $c$ with $0<c<1$. The function $\nu_\mathcal{P}:\ring{F}\rightarrow
\mathbb{R}^+$ given by $\nu_\mathcal{P}(x) = c^{n_\mathcal{P}(x)}$ is a
non-Archimedean valuation on $F$. Moreover, all non-Archimedean valuations on
$F$ are equivalent to a valuation of this form.

We summarize these facts in the following theorem.  \begin{thm} Let $F$ be a
number field. Then every Archimedean valuation of $F$ is equivalent to
$\nu_\sigma$ for some real or complex embedding $\sigma$ of $F$, and every
non-Archimedean valuation of $F$ is equivalent to $\nu_\mathcal{P}$ for some
prime ideal $\mathcal{P}$ of $\ring{F}$. The former are sometimes called
\emph{infinite places,} while the latter are called \emph{finite places}.
\end{thm}

A valuation $\nu$ on a field $F$ defines a metric on $F$ by $d(x,y) = \nu (x -
y)$. The completion of $F$ with respect to this metric is denoted by $F_\nu$.
Equivalent valuations give rise to the same completions. If $\nu = \nu_\sigma$
for a real or complex embedding $\sigma$ of $F$, then $F_\nu$ is isomorphic to
$\mathbb{R}$ or $\mathbb{C}$ respectively. If $\nu = \nu_\mathcal{P}$ for some
prime ideal $\mathcal{P}\subseteq \ring{F}$, then $F_\nu$ is called a
\emph{$\mathcal{P}$-adic field}.

Let $A$ be a quaternion algebra over a number field $F$, and let $F_\nu$ be the
completion of $F$ with respect to some valuation $\nu$. Then we can construct
the tensor product $A\otimes_F F_\nu$, which turns out to be a quaternion
algebra over $F_\nu$. Indeed, if $A = \hilb{a}{b}{F}$, then
$A\otimes_F F_\nu = \hilb{a}{b}{F_\nu}$. Ultimately, the
classification of quaternion algebras over $F$ will be reduced to the
classification of quaternion algebras over completions of $F$ with respect to
valuations $\nu$. The two following theorems will be useful in this regard.
\begin{thm} Let $\mathbb{R}$ be the real number field. Then $\mathcal{H}$ is
the unique quaternion division algebra over $\mathbb{R}$.  \end{thm}
\begin{thm} Let $F_\nu$ be a $\mathcal{P}$-adic field. Then there is a unique
quaternion division algebra over $F_\nu$.  \end{thm}

Thus when $A$ is a quaternion algebra over a number field $F$ and $A_\nu$ is
the corresponding quaternion algebra over a real or $\mathcal{P}$-adic
completion $F_\nu$ of $F$, by Theorem \ref{split} there are two possibilities:
$A$ is the unique quaternion division algebra over $F_\nu$, or $A \cong
M_2(F_\nu)$. In the former case, we say that $A$ \emph{ramifies} at $\nu$,
while in the latter case we say that $A$ \emph{splits} at $\nu$. When $F_\nu =
\mathbb{R}$, there is a simple test to determine which algebra is represented
by the Hilbert symbol $\hilb{a}{b}{\mathbb{R}}$:
\begin{equation}\label{realcondition}
\textrm{if $a$ and $b$ are both negative, then $A\cong\mathcal{H}$, otherwise $A$ splits.}
\end{equation}
Various tests exist for $\mathcal{P}_\nu$, but often the simplest test is to
determine if there exists a solution to equation (\ref{hilberteq}).  (See Appendix \ref{APP:RAMCHECK}.) Notice that
by Theorem \ref{THM:HILBERT_EQ} every quaternion algebra over the complex
numbers is isomorphic to $M_2(\mathbb{C})$.

The following theorem provides the necessary criterion for distinguishing between quaternion algebras over number fields.
\begin{thm}(Vign\'eras \cite{VIG}.)
Let $F$ be a number field. For each quaternion algebra $A$ over $F$, denote by $Ram(A)$ the set of all real or finite places at which $A$ ramifies. Then two quaternion algebras $A$ and $A'$ over $F$ are equal if and only if $Ram(A)=Ram(A')$.
\end{thm}

Thus the complete identification of a quaternion algebra $A$ over a number
field $F$ amounts to determining $Ram(A)$. It is easy to check if $A$ ramifies
at the real infinite places of $F$. Let $\alpha$ be a primitive element of $F$
(i.e. an element whose powers form a basis for $F$ over $\mathbb{Q}$.) Every
embedding of $F$ in $\mathbb{R}$ can be obtained from a real root $\alpha_i$ of
the minimal polynomial of $\alpha$ over $\mathbb{Q}$ by extending the map
$\sigma_i:\alpha \rightarrow \alpha_i$ linearly to $F$. If $A=\hilb{a}{b}{F}$
and $a$ and $b$ are expressed as polynomials in $\alpha$, then it is
straightforward to check if condition (\ref{realcondition}) holds for
$\sigma_i(a)$ and $\sigma_i(b)$.

The finite places of $F$ are more difficult to check.  The most straightforward
method is to check to see if the Hilbert equation (\ref{hilberteq}) has a
solution; our procedure for doing this is detailed in Appendix
\ref{APP:RAMCHECK}.

\subsection{Arithmeticity}

The notion of an arithmetic group comes from the theory of algebraic groups and
is a standard way of producing finite-covolume discrete subgroups of
semi-simple Lie groups.  To see how this general theory relates to Kleinian
groups, see \cite{HILDEN} or \cite{MR_BOOK}.

In the case of Kleinian groups, the following definition coincides with the
most general one, and is naturally related to the quaternion algebras which
we have already mentioned.

Let $A$ be a quaternion algebra over a number field $F$ and denote by $R_F$ the
ring of integers in $F$.  An {\em order} ${\cal O}$ in $A$ is an $R_F$-lattice
(spanning $A$ over $F$) that is also a ring with unity.  For every complex
place $\nu$ of $F$ there is an embedding of $A
\longrightarrow M_2(\mathbb{C})$ determined by the isomorphism $A \otimes_F
F_\nu \cong M_2(\mathbb{C})$.  Given a complex place $\nu$ and an order ${\cal
O}$, we can construct a subgroup of $\SL(2,\mathbb{C})$, and hence of
$\PSL(2,\mathbb{C})$, by taking the image $\Gamma_{\cal O}^\nu$ of the elements
of ${\cal O}$ with unit norm under the embedding $A \longrightarrow
M_2(\mathbb{C})$ defined above.

In the case that $F$ has a unique complex place $\nu$ and that $A$ ramifies over
every real place of $F$, then $\Gamma_{\cal O} := \Gamma_{\cal O}^\nu$ is a discrete
subgroup of $\PSL(2,\mathbb{C})$ (see Sections 8.1 and 8.2 of \cite{MR_BOOK}).

\begin{defn}
A Kleinian group $\Gamma$ is called \emph{arithmetic} if it is commensurable
with $\Gamma_{\cal O}$ for some order ${\cal O}$ of a quaternion algebra
that ramifies over every real place and is defined over a field with a unique
complex place.
\end{defn}

Viewing $\Isom^+(\mathbb{H}^3)$ as $SO^+(3,1)$ furnishes an alternative
construction of arithmetic Kleinian groups as follows. Let $F$ be a real number
field, and $(V,q)$ a four-dimensional quadratic space over $F$ with signature
$(3,1)$. Any $F$-linear map $\sigma : V \longrightarrow V$ preserving $q$ can
be identified with an element of $SO^+(3,1)$, and thus of
$\Isom^+(\mathbb{H}^3)$, by extension of scalars from $F$ to $\mathbb{R}$.
Given an $\ring{F}$-lattice $L \subset V$ of rank 4 over $\ring{F}$, the group
$\SO(L) := \{ \sigma \in SO^+(3,1) \cap GL(4,F) : \sigma(L) = L \}$ is always
discrete and arithmetic.

Moreover, the groups of the form $\SO(L)$ give representatives for the
commensurability classes of all Kleinian groups that possess a non-elementary
Fuchsian subgroup. Since all reflection groups determined by finite-volume
polyhedra have non-elementary Fuchsian subgroups, for our purposes this can be
considered the definition of arithmeticity. See \cite[page 143]{HILDEN} for a
discussion. The distinction between arithmetic groups arising from quaternion
algebras and those arising from quadratic forms is also discussed in
\cite[pages 217--221]{VS}, whose authors call the latter ``arithmetic groups of
the simplest kind''.

Aside from its relationship to algebraic groups, arithmeticity is interesting
for many reasons including the fact that for arithmetic groups $\Gamma$ the
pair $(k \Gamma, A \Gamma)$ is a complete commensurability invariant (see
Section 8.4 of \cite{MR_BOOK}). This will allow us to identify several
unexpected pairs of commensurable reflection groups which are presented in Section
\ref{SUBSEC:COMMENS_PAIRS}.

To decide whether a given reflection group $\Gamma$ determined by a polyhedron
$P$ is arithmetic there is a classical theorem due to Vinberg \cite{VIN}:

\begin{thm}\label{THM:ARITH}
  Let $(a_{ij})$ be the Gram matrix of a Coxeter polyhedron $P$.
  Then the reflection group determined by $P$ is arithmetic if and only if
  the following three conditions hold:
  \begin{enumerate}
    \item $K := \mathbb{Q}(a_{ij})$ is totally real.
    \item For every embedding $\sigma : K \longrightarrow \mathbb{C}$ such that $\sigma |_{k(P)} \neq \id$ (where $k(P)$ is the field defined in Theorem \ref{THM:POLY_ITF}), the matrix $(\sigma(a_{ij}))$ is positive semi-definite.
    \item The $a_{ij}$ are algebraic integers.
  \end{enumerate}
\end{thm}

More generally, for any finite-covolume Kleinian group $\Gamma$, Maclachlan and Reid have proved a similar result,
Theorem 8.3.2 in \cite{MR_BOOK}:

\begin{thm}\label{THM:ARITH2}
  A finite-covolume Kleinian group $\Gamma$ is arithmetic if and only if the following three conditions hold:
  \begin{enumerate}
    \item $k \Gamma$ has exactly one complex place.
    \item $A \Gamma$ ramifies at every real place of $k\Gamma$.
    \item $tr \gamma$ is an algebraic integer for each $\gamma \in \Gamma$.
  \end{enumerate}
\end{thm}


\section{Worked example}\label{SEC:WORKED_EXAMPLE}

A Lambert cube is a compact polyhedron realizing the combinatorial type of a
cube, with three disjoint non-coplanar edges chosen and assigned dihedral angles
$\frac{\pi}{l},$ $\frac{\pi}{m}$, and $\frac{\pi}{n}$, and the remaining edges
assigned dihedral angles $\frac{\pi}{2}$.  It is easy to verify that if $l,
m, n > 2$, then, such an assignment of dihedral angles  satisfies the
hypotheses of Andreev's Theorem.  The resulting polyhedron is called the
$(l,m,n)$-Lambert Cube, which we will denote by
$P_{l,m,n}$.  

In this section we illustrate our techniques by computing the invariant trace
field and invariant quaternion algebra associated to the $(3,3,6)$ Lambert
cube.

The starting point of our computation is a set of low-precision decimal
approximations of outward-pointing normal vectors $\{{\bf e}_1,\ldots, {\bf
e}_6\}$ to the six faces of our cube.  For a given compact hyperbolic
polyhedron, it is nontrivial to construct such a set of outward-pointing normal
vectors.  One way is to use the collection of Matlab scripts described in
\cite{ROE_CONSTRUCTION}.  Throughout this paper we will always assume the
following normalization for the location of our polyhedron: the first three faces meet at a vertex, the first face
has normal vector $(0,0,0,*)$, the second face has normal vector $(0,0,*,*)$,
and the third has form $(0,*,*,*)$, where $*$ indicates that no condition is
placed on that number. 

We then use Newton's Method with extended-precision decimals to improve this set
of approximate normals until they are very precise.  (Here, we do this with
precision $40$ numbers, but we display fewer digits for the reader.)  The vectors
$\{{\bf e}_1,\ldots {\bf e}_6\}$ are displayed as rows in the following matrix:

{\small
\begin{eqnarray*}
\left[ \begin {array}{cccc}  0.0&- 0.0& 0.0&-
0.99999996237702\\\noalign{\medskip} 0.0& 0.0& 0.86602540463740&
0.50000001881149
\\\noalign{\medskip}- 0.0&- 1.00000002377892&- 0.0&{ 1.0\times 10^{-14
}}\\\noalign{\medskip} 1.38941010745090& 0.86602538319131&-
 0.73831913868376& 1.27880621354777\\\noalign{\medskip}
 0.79708547435960& 1.27880618290479& 0.0&- 0.0\\\noalign{\medskip}
 0.62728529885922& 0.0&- 1.18046043888844&{ 1.0\times 10^{-14}}
\end {array} \right]
\end{eqnarray*}
}

The normalization we have chosen for the location of our polyhedron assures us
that each of these decimals should approximate an algebraic number of some
(low) degree.  There are commands in many computer algebra packages for
guessing the minimal polynomial that is most likely satisfied by a given
decimal approximate.  Most of these commands are ultimately based on the LLL
algorithm \cite{LLL}.  (We have used the command {\tt minpoly()} in Maple, the
command {\tt algdep()} in Pari/GP, and the command {\tt RootApproximant[]} in
Mathematica 6.) Each of these commands requires a parameter specifying up to
what degree of polynomials to search.  In this case we specify degree $30$.
The resulting matrix of guessed minimal polynomials is:

{\small
\begin{eqnarray*}
\left[ \begin {array}{cccc} {\it X}&{\it X}&{\it X}&1+{\it
X}\\\noalign{\medskip}{\it X}&{\it X}&-3+4\,{{\it
X}}^{2}&-1+2\,{
\it X}\\\noalign{\medskip}{\it X}&1+{\it X}&{\it X}&{\it X}
\\\noalign{\medskip}-9-88\,{{\it X}}^{2}+48\,{{\it X}}^{4}&-3+4\,{
{\it X}}^{2}&1-28\,{{\it X}}^{2}+48\,{{\it X}}^{4}&3-28\,{{\it
X}}^{2}+16\,{{\it X}}^{4}\\\noalign{\medskip}-9+4\,{{\it X}}^{2}
+16\,{{\it X}}^{4}&3-28\,{{\it X}}^{2}+16\,{{\it X}}^{4}&{\it
X}&{\it X}\\\noalign{\medskip}-1-{{\it X}}^{2}+9\,{{\it X}}^{4
}&{\it X}&-3+{\it X}+3\,{{\it X}}^{2}&{\it X}\end {array}
 \right]
\end{eqnarray*}
}

The next step is to specify which root of each given minimal polynomial is
closest to the decimal approximate above.  In the current case, the solutions
of each polynomial are easily expressed by radicals, so we merely pick the
appropriate expression.  (For more complicated examples, our computer program
uses a more sophisticated way of expressing algebraic numbers as described 
in Section \ref{SEC:PROGRAM}.)  For the current example, the following 
matrix contains as rows our guessed exact values for 
$\{{\bf e}_1,\ldots {\bf e}_6\}$.

{\small
\begin{eqnarray*}
N:=
\left[ \begin {array}{cccc} 0&0&0&-1\\\noalign{\medskip}0&0&\frac{\sqrt
{3}}{2}&\frac{1}{2}\\\noalign{\medskip}0&-1&0&0\\\noalign{\medskip}
\frac{\sqrt {33+6\,\sqrt {37}}}{6}&\frac{\sqrt {3}}{2}&\frac{-\sqrt {42+6\,\sqrt {37}
}}{12}&\frac{\sqrt {14+2\,\sqrt {37}}}{4}\\\noalign{\medskip}\frac{\sqrt {-2+2\,
\sqrt {37}}}{4}&\frac{\sqrt {14+2\,\sqrt {37}}}{4}&0&0\\\noalign{\medskip}
\frac{\sqrt {2+2\,\sqrt {37}}}{6}&0&\frac{-1-\sqrt {37}}{6}&0\end {array} \right]
\end{eqnarray*}
}

Corresponding to this set of guessed normal vectors we have the Gram Matrix $G_{i,j} = 2\langle {\bf e}_i,{\bf e}_j\rangle$:
{\small
\begin{eqnarray*} G:=
 \left[ \begin {array}{cccccc} 2&-1&0&\frac{-\sqrt {14+2\,\sqrt {37}}}{2}&0
&0\\\noalign{\medskip}-1&2&0&0&0&\frac{-\sqrt {3} \left( 1+\sqrt {37}
 \right)}{6} \\\noalign{\medskip}0&0&2&-\sqrt {3}&\frac{-\sqrt {14+2\,
\sqrt {37}}}{2}&0\\\noalign{\medskip}\frac{-\sqrt {14+2\,\sqrt {37}}}{2}&0&-
\sqrt {3}&2&0&0\\\noalign{\medskip}0&0&\frac{-\sqrt {14+2\,\sqrt {37}}}{2}&0
&2&-1\\\noalign{\medskip}0&\frac{-\sqrt {3} \left( 1+\sqrt {37}
 \right)}{6} &0&0&-1&2\end {array} \right]
\end{eqnarray*}
}

By checking that there are $2$'s down the diagonal of $G$ and that there is
$-2\cos(\alpha_{ij})$ in the $ij$-th entry of $G$ if faces $i$ and $j$ are
adjacent, we can see that the guessed matrix $N$ whose rows represent
outward-pointing normal vectors was correct.  Here there are $4\cdot6-6 = 18$
equations that we have checked, consistent with the number of guessed values in
the matrix $N$.  Consequently $G$ is the exact Gram matrix for the $(3,3,6)$
Lambert cube.  We can use this to compute the invariant trace field and the
invariant quaternion algebra.

\begin{figure}
\begin{center}
\begin{picture}(0,0)%
\epsfig{file=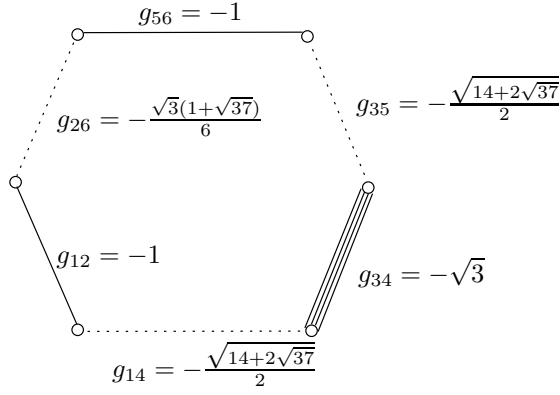}%
\end{picture}%
\setlength{\unitlength}{3947sp}%
\begingroup\makeatletter\ifx\SetFigFont\undefined%
\gdef\SetFigFont#1#2#3#4#5{%
  \reset@font\fontsize{#1}{#2pt}%
  \fontfamily{#3}\fontseries{#4}\fontshape{#5}%
  \selectfont}%
\fi\endgroup%
\begin{picture}(2303,2390)(3326,-3858)
\put(5510,-2145){\makebox(0,0)[lb]{\smash{{\SetFigFont{10}{12.0}{\familydefault}{\mddefault}{\updefault}{\color[rgb]{0,0,0}$g_{35}=-\frac{\sqrt{14+2\sqrt{37}}}{2}$}%
}}}}
\put(3628,-3098){\makebox(0,0)[lb]{\smash{{\SetFigFont{10}{12.0}{\familydefault}{\mddefault}{\updefault}{\color[rgb]{0,0,0}$g_{12}=-1$}%
}}}}
\put(4141,-1569){\makebox(0,0)[lb]{\smash{{\SetFigFont{10}{12.0}{\familydefault}{\mddefault}{\updefault}{\color[rgb]{0,0,0}$g_{56}=-1$}%
}}}}
\put(3975,-3815){\makebox(0,0)[lb]{\smash{{\SetFigFont{10}{12.0}{\familydefault}{\mddefault}{\updefault}{\color[rgb]{0,0,0}$g_{14}=-\frac{\sqrt{14+2\sqrt{37}}}{2}$}%
}}}}
\put(3629,-2241){\makebox(0,0)[lb]{\smash{{\SetFigFont{10}{12.0}{\familydefault}{\mddefault}{\updefault}{\color[rgb]{0,0,0}$g_{26}=-\frac{\sqrt{3}(1+\sqrt{37})}{6}$}%
}}}}
\put(5522,-3215){\makebox(0,0)[lb]{\smash{{\SetFigFont{10}{12.0}{\familydefault}{\mddefault}{\updefault}{\color[rgb]{0,0,0}$g_{34}=-\sqrt{3}$}%
}}}}
\end{picture}%

\end{center}
\caption{\label{FIG:COX} Coxeter symbol for the $(3,3,6)$ Lambert cube.  Next to each edge we display the corresponding element of the Gram matrix $G$.}
\end{figure}

The nontrivial cyclic products from $G$ correspond to non-trivial cycles in the Coxeter symbol for
$P_{3,3,6}$, which is depicted in Figure \ref{FIG:COX} with the appropriate element of the Gram matrix written next to each 
edge.  (For those unfamiliar with Coxeter symbols, see \ref{SEC:PROGRAM}.)

The non-trivial cyclic products correspond to closed loops in the Coxeter symbol.  Always included are  the squares of each
entry of $G$, of which the only two irrational ones are: $(g_{12})^2 =
\frac{7+\sqrt{37}}{2}$ and $(g_{35})^2 = \frac{19+\sqrt{37}}{6}$.
The other nontrivial cyclic product corresponds to the closed loop in
the Coxeter symbol: $g_{12}g_{26}g_{65}g_{53}g_{34}g_{41} = 11+2\sqrt{37}$. 

Thus, $k(P_{3,3,6}) = \mathbb{Q}(\sqrt{37})$.

Notice that ${\bf v}_1 := {\bf e}_1$, ${\bf v}_{12} := g_{12}{\bf e}_2 = -{\bf
e}_2$, ${\bf v}_{143} := g_{14}g_{43}{\bf e}_3 = \frac{-\sqrt {14+2\,\sqrt
{37}}}{2}\cdot \sqrt{3} {\bf e}_3$, and ${\bf v}_{14} := g_{14}{\bf e}_4 = \frac{-\sqrt {14+2\,\sqrt
{37}}}{2} {\bf e}_4$ are linearly independent, so that they span the quadratic space
$(M,q)$.  (See the description of Theorems \ref{THM:POLY_ITF} and \ref{THM:POLY_QA}.) We now
compute the matrix representing $q$ with respect to the basis
$\{{\bf v}_1,{\bf v}_{12},{\bf v}_{143},{\bf v}_{14}\}$.

{\small
\begin{eqnarray*}
\left[ \begin {array}{cccc} 2&1&0&\frac{7+\sqrt {37}}{2}
\\\noalign{\medskip}1&2&0&0\\\noalign{\medskip}0&0&21+3\,\sqrt {37}&\frac{21
+3\sqrt {37}}{2}\\\noalign{\medskip}\frac{7+\sqrt {37}}{2}&0&\frac{21+3
\sqrt {37}}{2}&7+\sqrt {37}\end {array} \right] 
\end{eqnarray*}
}

The determinant of this matrix, and hence the discriminant of the quadratic
form $q$ (a number well-defined up to a square in the field $k(P_{3,3,6})$) is
$d=\frac{2973-489\sqrt{37}}{2}$, consequently $k\Gamma_{3,3,6} =
k(P_{3,3,6})\left(\sqrt{ \frac{2973-489\sqrt{37}}{2}}\right)$.  Since $d$ is
primitive for $k(P_{3,3,6})$, we actually have $k\Gamma_{3,3,6} =
\mathbb{Q}\left(\sqrt{\frac{2973-489\sqrt{37}}{2}}\right)$.  This expression still looks
rather cumbersome, and by writing a minimal polynomial for
$\sqrt{\frac{2973-489\sqrt{37}}{2}}$ and using the ``polredabs()'' command in
Pari, we can check that $\sqrt{-10-2\,\sqrt {37}}$  also generates this field,
hence $k\Gamma_{3,3,6} = \mathbb{Q}\left(\sqrt{-10-2\,\sqrt {37}}\right)$.

\vspace{0.1in}
In order to use Theorems \ref{THM:POLY_ITF} and \ref{THM:POLY_QA} to compute the invariant quaternion algebra we need to express the quadratic form $q$ with respect to
an orthogonal basis $\{{\bf w}_1, {\bf w}_2, {\bf w}_3, {\bf w}_4\}$.  The result is:

{\small
\begin{eqnarray*}
 \left[ \begin {array}{cccc} 2&0&0&0\\\noalign{\medskip}0&\frac{3}{2}&0&0
\\\noalign{\medskip}0&0&21+3\,\sqrt {37}&0\\\noalign{\medskip}0&0&0&
\frac {-151-25\sqrt {37}}{12}\end {array} \right]
\end{eqnarray*}
}

Thus, a Hilbert symbol describing $A\Gamma_{3,3,6}$ is given by $\hilb{-q({\bf w}_1)q({\bf w}_2)}{
-q({\bf w}_1)q({\bf w}_3)}{k\Gamma_{3,3,6}}$.
That is, $A\Gamma_{3,3,6} \cong
\hilb{-3}{-42-6\sqrt{37}}{\mathbb{Q}\left(\sqrt{-10-2\,\sqrt {37}}\right)}$.

Without using the machinery described in Section \ref{SEC:QUATERNION} is
difficult to interpret this Hilbert symbol.  The invariant
trace field $k\Gamma_{3,3,6}$ has two real places, and $A\Gamma_{3,3,6}$ is
ramified at each of these two places.  By attempting to solve the Hilbert
Equation (see Appendix \ref{APP:RAMCHECK}) we also observe that $A\Gamma_{3,3,6}$ is
ramified over exactly two finite prime ideals in the ring of integers from
$k\Gamma_{3,3,6}$.  These prime ideals lie over the rational prime $3$ and we
denote them by ${\cal P}_3$ and ${\cal P'}_3$.  Thus,
according to Vign\'eras \cite{VIG}, this finite collection of ramification data provides a
complete invariant for the isomorphism class of  $A\Gamma_{3,3,6} \cong
\hilb{-3}{-42-6\sqrt{37}}{\mathbb{Q}\left(\sqrt{-10-2\,\sqrt {37}}\right)}$.

Using either Theorem \ref{THM:ARITH} or Theorem \ref{THM:ARITH2}, one can also
deduce that $\Gamma_{3,3,6}$ is not arithmetic because the Gram matrix contains
the element $\frac{-\sqrt {14+2\,\sqrt {37}}}{2}$, which is not an algebraic
integer.  However the other two hypotheses of each of these theorems are
satisfied.

\section{Description of the program}\label{SEC:PROGRAM}

We've written a collection of PARI/GP scripts that automate the procedure from Section
\ref{SEC:WORKED_EXAMPLE} for (finite-volume) hyperbolic polyhedra. The scripts take
as input two matrices: the matrix of face normals, whose rows are
low-precision decimal approximations to the outward-pointing normal vectors
(normalized in $\mathbb{H}^3$ as in the previous section); and the matrix of edge labels, a
square matrix $(n_{ij})$ whose diagonal entries must be one, and whose
off-diagonal terms $n_{ij}$ describe the relation between the $i$-th and $j$-th
faces: $n_{ij} = 0$ means the faces are non-adjacent and any other value means
that they meet at a dihedral angle of $\pi/n_{ij}$. The matrix of edge labels
actually determines the polyhedron uniquely (up to hyperbolic isometry); one can
use \cite{ROE_CONSTRUCTION} to compute the approximate face normals from the
edge labels.

In a typical session, after loading the input matrices, the user runs Newton's
Method to obtain a higher-precision approximation for the face normals, and
then can have the computer guess and verify exact values for the normal vectors
and for the Gram matrix. With the exact Gram matrix, the user can request
the invariant trace field (described by a primitive element) and the invariant
quaternion algebra (described by a Hilbert symbol, or, after an additional
command, by ramification data). There is also a function to test the polyhedron
for arithmeticity---this too takes the Gram matrix as input.

The scripts are available at \cite{SNAPHEDRON}; the package includes a sample session
and a user guide. The scripts, as mentioned before, are written in the high-level
language GP, which helped in our effort to make the source code as readable as possible.
In fact, the reader can consider the source code as executable statements of the theorems
in the previous sections.

There are two main technical challenges involved in writing these scripts:
choosing an appropriate representation of general algebraic numbers, and
choosing a systematic
way of listing the non-zero cyclic products $a_{i_1 i_2} a_{i_2 i_3} \cdots
a_{i_r i_1}$ from the Gram matrix.

We describe a given algebraic number $\alpha$ as a pair $(p(z),\tilde
\alpha)$, where $p(z)$ is the minimal polynomial for $\alpha$ over
$\mathbb{Q}$ and $\tilde \alpha$ is a decimal approximation to
$\alpha$.  This representation is rather common---it is described in
the textbook on computational number theory by Cohen \cite{COHEN}.
While this representation is not already available in PARI/GP
\cite{PARI2}, it is easy to program an algebraic number package
working within PARI/GP for this representation.

A disadvantage of this representation for algebraic numbers is that
when performing arithmetic on algebraic numbers, it is usually
necessary to find composita of the fields generated by each number.
This is not only a programming difficulty, it is the slowest part of
our program.  In some cases, when we know that we will do arithmetic
with a given set of numbers $\{\alpha_1,\ldots,\alpha_n\}$, we are able
to speed up our calculations by computing a primitive element $\beta$
for $\mathbb{Q}\left(\alpha_1,\ldots,\alpha_n\right)$ and
re-expressing each $\alpha_i$ as an element of $\mathbb{Q}[z]/\langle
p(z) \rangle$ where $p(z)$ is the minimal polynomial for $\beta$ over
$\mathbb{Q}$.  (PARI/GP has a data type called ``polmod'' for this
representation.) Once expressed in terms of a common field, algebraic
computations in terms of these polmods are extremely fast in PARI/GP.

Of course we don't typically have {\em a~priori} knowledge of the
outward-pointing normal vectors or of the entries in the Gram matrix
for our polyhedron $P$.  However, decimal approximations can be
obtained using \cite{ROE_CONSTRUCTION}.  Just like in Section
\ref{SEC:WORKED_EXAMPLE}, we use the LLL algorithm \cite{LLL} to guess
minimal polynomials for the algebraic numbers represented by these
decimal approximations.  Typically a rather high-precision
approximation is needed (sometimes $100$ digits of precision) to
obtain a correct guess.  In this case, we use Newton's Method and the high-precision
capabilities of PARI/GP to improve the precision
of the normal vectors obtained from \cite{ROE_CONSTRUCTION}.
Correctness of the guess is verified once we use the guessed algebraic
numbers for the outward-pointing normals to compute the Gram matrix
and verify that each entry corresponding to a dihedral angle
$\frac{\pi}{n}$ has the correct minimal polynomial for
$-2\cos\left(\frac{\pi}{n}\right)$.  Additionally we check
that the diagonal entries are exactly $2$.  In some
cases the guessed polynomials are not correct, but typically, by
sufficiently increasing the number of digits of precision, reapplying
Newton's Method, and guessing again, we arrive at correct guesses.

Note that the number of checks, i.e. equations, equals the number of variables.
If $P$ has $n$ faces, there are $4n-6$ variables: one for each coordinate of
each outward-pointing normal, minus the six coordinates normalized to $0$.
Since we deal with compact Coxeter polyhedra and these necessarily have three
faces meeting at each vertex, the number of edges is $3n-6$; there is one
equation for each of these, and there are $n$ additional equations for the
diagonal entries of the Gram matrix, giving a total of $4n-6$ equations.

The guess and check philosophy is inspired by SNAP \cite{SNAP}, which
also goes through the process of improving an initial approximation to
the hyperbolic structure, guessing minimal polynomials, and verifying.
Such techniques are also central to many areas of experimental
mathematics in which the LLL algorithm is used to guess linear
dependencies that are verified {\em a posteriori}.

When it comes to listing the non-zero cyclic products from the Gram matrix, we
try to avoid finding more of them than necessary to compute the field they
generate. It is easier to discuss these products in terms of the Coxeter symbol
for the polyhedron.  Recall that the Coxeter symbol is a graph whose vertices
are the faces of the polyhedron with edges between pairs of non-adjacent faces
and also between pairs of adjacent non-perpedicular faces. Typically edges
between non-adjacent faces are drawn dashed, and edges between adjacent ones
meeting at an angle of $\pi/n$ are labeled $n-2$.   A non-zero cyclic product in the Gram matrix
corresponds to a closed path in the Coxeter symbol. 

An example Coxeter symbol appears in Figure \ref{FIG:COX}.  While this Coxeter
symbol and most that appear in the literature are planar, this is not typically
the case due to the large number of dashed edges.  The reason so many planar
Coxeter symbols appear in the literature is that they are especially useful for
tetrahedra in high dimensions and those are planar.

Our method, then, is to list all the squares of the elements of the Gram
matrix, one for each edge in the Coxeter symbol, and a set of cycles that form
a basis for the $\mathbb{Z} / 2\mathbb{Z}$-homology of the Coxeter symbol.
These allow one to express any cyclic product in the Gram matrix as a product
of a number of basic cycles and squares or inverses of squares of Gram matrix
elements corresponding to edges traversed more than once. 

To get such a basis one can take any spanning tree for the graph and then, for
each non-tree edge, take the cycle formed by that edge and the unique tree-path
connecting its endpoints.  Finding a spanning tree for a graph is a classical
computer science problem for which we use breadth-first search (see
\cite{COMPSCI}). The spanning tree is also used for the task of finding a basis
for $M(P)$, as this requires finding paths in the Coxeter symbol from a fixed
vertex to four others. Breadth-first search has the advantage of producing a
short bushy tree which in turn gives short paths and cycles.

\section{Many computed examples}\label{SEC:RESULTS}
In this section we present a number of examples including some unexpectedly commensurable pairs of groups.   We also present some borderline cases of groups that are incommensurable, but have
some of the invariants in common.

\subsection{Lambert Cubes}

Recall from Section \ref{SEC:WORKED_EXAMPLE} that a Lambert cube is
a compact polyhedron realizing the combinatorial type of a
cube, with three disjoint non-coplanar edges chosen and assigned dihedral angles
$\frac{\pi}{l},$ $\frac{\pi}{m}$, and $\frac{\pi}{n}$ with $l,m,n > 2$, and the remaining edges
assigned dihedral angles $\frac{\pi}{2}$.  Any reordering of $(l,m,n)$ can be obtained by applying 
an appropriate (possibly orientation reversing) isometry, so when studying Lambert cubes it suffices to
consider triples with $l \leq m \leq n$.

In Table \ref{TAB:LAMBERT} we provide the invariant trace fields and the
ramification data for the invariant quaternion algebras for Lambert cubes with
small $l,m,$ and $n$.

\begin{table}
\begin{center}
\scalebox{0.7}{
\begin{tabular}{|l|l|l|l|l|l|}
\hline
$(l,m,n)$ & $k\Gamma$ & disc & $A\Gamma$ ramification &  Arith? \\
\hline
\hline
$(3,3,3)$ &  $x^4 - x^3 - x^2 - x + 1$ & $-507$ & $[1, 1, -1, -1]$ & Yes \\
\hbox{}	 & $-0.651387818 - 0.758744956i$ & \hbox{} & $\emptyset$ & \hbox{} \\
\hline
$(3,3,4)$ & $x^2 + 1$ & $-4$ & $[-1, -1]$ & No*  \\
\hbox{}	 & $1.00000000000000i$ & \hbox{} & $\emptyset$ &  \\
\hline
$(3,3,5)$ & $x^8 - x^7 - 3x^6 + x^4 - 3x^2 - x + 1$ &   $102378125$ & $[1, 1, 1, 1, -1, -1, -1, -1]$ & No\\
\hbox{} & $0.725191949 - 0.688546757i$  & \hbox{} &  $\emptyset$ & \hbox{} \\
\hline
$(3,3,6)$ & $x^4 + 5x^2 - 3$  &  $-4107$ & $[1, 1, -1, -1]$ & No* \\
\hbox{} & $- 2.354013863i$ & \hbox{} & ${\cal P}_3, {\cal P}_3'$ &  \\
\hline
$(3,4,4)$ & $x^4 - 2x^3 - 2x + 1$ & $-1728$ & $[1, 1, -1, -1]$ & Yes \\
\hbox{} & $-0.3660254039 + 0.930604859i$ & \hbox{} & $\emptyset$ & \\
\hline
$(3,4,5)$ & $x^8 - x^7 - 18x^6 - 18x^5 + 95x^4 + 218x^3 + 182x^2 + 71x + 11$  & $249761250000$ & $[1, 1, 1, 1, -1, -1, -1, -1]$ & No \\
\hbox{} & $-2.142201597 + 1.146040793i$ & \hbox{} & $\emptyset$ & \\
\hline
$(3,4,6)$  & $x^4 - x^3 - 11x^2 + 33x - 6$ & $-191844$ &  $[1, 1, -1, -1]$ & No*\\
& $2.386000936 - 1.441874268i$ &  & ${\cal P}_2, {\cal P}_2'$ &  \\
\hline

$(3,5,5)$ & $x^8 - 3x^7 + x^5 + 3x^4 + x^3 - 3x + 1$ & $184280625$ & $[1, 1, 1, 1, -1, -1, -1, -1]$ & No \\
\hbox{} & $-0.2251919494 + 0.974314418i$ & \hbox{} & $\emptyset$ & \\
\hline
$(3,5,6)$ & $x^8 - x^7 - 19x^6 + 30x^5 + 74x^4 - 72x^3 - 310x^2 + 413x - 71$  & $876653128125$ & $[1, 1, 1, 1, -1, -1, -1, -1]$  & No  \\
\hbox{} & $-1.6937824749171233274 - 1.2086352565016507434i$ &  & $\emptyset$ &  \\
\hline
$(3,6,6)$ &  $x^4 - x^3 - x^2 - x + 1$ & $-507$ & $[1, 1, -1, -1]$ & Yes \\
\hbox{}	 & $-0.651387818 - 0.758744956i$ & \hbox{} & $\emptyset$ & \hbox{} \\
\hline
$(4,4,4)$ & $x^4 - x^2 - 1$ & $-400$ & $[1, 1, -1, -1]$ & Yes \\
\hbox{} & $0.786151377i$ & \hbox{} & $\emptyset$ & \\
\hline
$(4,4,5)$ & $x^8 - 4x^7 + 4x^6 - 6x^5 + 19x^4 - 14x^3 + 4x^2 - 6x + 1$ & $1548800000$ & $[1, 1, 1, 1, -1, -1, -1, -1]$ & No \\
\hbox{} & $-0.748606020 + 1.434441708i$ & \hbox{} & $\emptyset$ & \\
\hline
$(4,4,6)$ & $x^4 - 2x^3 - 6x + 9$ &  $-9408$ & $[1, 1, -1, -1]$ & No* \\
\hbox{} & $-0.822875655 + 1.524098309i$ & \hbox{} & ${\cal P}_3, {\cal P}_3'$ &  \\
\hline
$(4,5,5)$ & $x^8 + 5x^6 - 3x^4 - 20x^2 + 16$ &  $2714410000$ & $[1, 1, 1, 1, -1, -1, -1, -1]$ & No\\
\hbox{} & $1.675405432i$ & \hbox{} & ${\cal P}_2, {\cal P}_2'$ & \\
\hline                
$(4,5,6)$ & $x^8 - 3x^7 - 25x^6 + 95x^5 + 50x^4 - 3x^3 - 1751x^2 + 2600x - 995$ & $21059676450000$ & $[1, 1, 1, 1, -1, -1, -1, -1]$ & No \\
\hbox{} & $-1.713666 - 1.949254i$ & \hbox{} & $\emptyset$ & \\
\hline
$(4,6,6)$ & $x^4 - x^3 - x^2 + 7x - 2$ & $-10404$ & $[1, 1, -1, -1]$ & No*\\
\hbox{} & $1.280776406 - 1.386060824i$ & \hbox{} & ${\cal P}_2, {\cal P}_2'$ &  \\
\hline
$(5,5,5)$ & $x^4 - x^3 + x^2 - x + 1$ & $125$ & $[-1, -1, -1, -1]$ & No \\
\hbox{} & $0.809016994 - 0.587785252i$ & \hbox{} & $\emptyset$ & \\
\hline
$(5,5,6)$ & $x^8 + 9x^6 + 13x^4 - 27x^2 + 9$ & $5863730625$ & $[1, 1, 1, 1, -1, -1, -1, -1]$ & No \\
\hbox{} & $-2.027709945i$ & \hbox{} & ${\cal P}_3, {\cal P}_3',{\cal P}_5, {\cal P}_5'$ & \\
\hline
$(5,6,6)$ & $x^8 - 2x^7 - 2x^6 + 16x^5 - 15x^4 - 19x^3 + 43x^2 - 22x + 1$ & $10637578125$ & $[1, 1, 1, 1, -1, -1, -1, -1]$ & No\\
\hbox{} & $1.353181081 - 1.579017768i$ & \hbox{} & $\emptyset$ & \\
\hline
$(6,6,6)$ & $x^4 - x^3 - 3x^2 - x + 1$ & $-1323$ & $[1, 1, -1, -1]$ & Yes \\
\hbox{} & $-0.895643923 - 0.4447718088i$ & \hbox{} & $\emptyset$ & \\
\hline
\end{tabular}
}
\end{center}
\caption{\label{TAB:LAMBERT} Commensurability invariants for Lambert Cubes.}
\end{table}

It is interesting to notice that many of the above computations are done by
hand in \cite{HILDEN}, where the authors determine which of the ``Borromean
Orbifolds'' are arithmetic by recognizing that they are $8$-fold covers of
appropriate Lambert cubes.  (They call the Lambert cubes ``pyritohedra.'')
Our results are consistent with theirs.

In this table and all that follow, we specify the invariant trace field
$k\Gamma$ by a canonical minimal polynomial $p(z)$ for a primitive element of
the field and, below this polynomial, a decimal approximation of the root that
corresponds to this primitive element.  

We specify the real ramification data for the invariant quaternion algebra
$A\Gamma$ by a vector ${\bf v}$ whose length is the degree of $p(z)$.  If the
$i$-th root of $p(z)$ (with respect to Pari's internal root numbering scheme)
is real, we place in the $i$-coordinate of ${\bf v}$ a $1$ if $A\Gamma$ is
ramified over corresponding completion of $k\Gamma$, otherwise we place a $0$.
If the $i$-th root of $p(z)$ is not real then we place a $-1$ in the
$i$-coordinate of ${\bf v}$.

Below this vector ${\bf v}$ we provide an indication of whether $A\Gamma$
ramifies at a prime ideal ${\cal P}_p$ over a rational prime $p$.  In the case
that $A\Gamma$ ramifies at multiple prime ideals over the same $p$, we list
multiple symbols ${\cal P}_p, {\cal P}'_p,\ldots$, etc.  More detailed
information about the generators of these prime ideals can be obtained in
PARI's internal format using our scripts.

The column labeled ``arith?'' indicates whether the group is arithmetic.  In
the case that the group satisfies conditions (1) and (2) from Theorems
\ref{THM:ARITH} and \ref{THM:ARITH2} but has elements with non-integral traces,
we place a star next to the indication that the group is not arithmetic.  (Some
authors refer to such groups as ``psuedo-arithmetic.'')

\subsection{Truncated cubes}

Another very simple family of hyperbolic reflection groups is obtained by
truncating a single vertex of a cube, assigning dihedral angles $\frac{\pi}{l},
\frac{\pi}{m},$ and $\frac{\pi}{n}$ to the edges entering the vertex that was
truncated and $\frac{\pi}{2}$ dihedral angles at all of the remaining edges.
Since the three edges entering the vertex that was truncated form a {\em
prismatic 3-circuit}, Andreev's Theorem provides the necessary and sufficient
condition that $\frac{1}{l} + \frac{1}{m} + \frac{1}{n} < 1$ for the existence
of such a polyhedron.  In Table \ref{TAB:TRUNCATED_CUBES} we list the $(l,m,n)$
truncated cubes for $l,m,n \leq 6$.

\begin{table}
\begin{center}
\scalebox{0.7}{
\begin{tabular}{|l|l|l|l|l|l|}
\hline
$(l,m,n)$ & $k\Gamma$ & disc & $A\Gamma$ ramification &  Arith? \\
\hline
\hline
$(3,3,4)$ & $x^4 - 2x^3 + x^2 + 2x - 1$ & $-448$ & $[1, 1, -1, -1]$ & Yes  \\
\hbox{}  & $1.207106781 - 0.978318343i$ & \hbox{} & $\emptyset$ &  \\
\hline
$(3,3,5)$ & $x^4 - x^2 - 1$ &   $-400$ & $[1, 1, -1, -1]$ & Yes\\
\hbox{} & $0.786151377i$  & \hbox{} &  $\emptyset$ & \hbox{} \\
\hline
$(3,3,6)$ & $x^4 + 2x^2 - 11$ & $-6336$   & $[1, 1, -1, -1]$ & No* \\
\hbox{} & $ 2.112842071i$  &  & $\emptyset$ &   \\
\hline
$(3,4,4)$ & $x^2 + 2$  &  $-8$ & $[-1, -1]$ & Yes \\
\hbox{} & $1.414213562i$ & \hbox{} & ${\cal P}_3, {\cal P}_3'$ & \\
\hline
$(3,4,5)$ & $x^8 + 14x^6 + 57x^4 + 86x^2 + 41$  & $26869760000$  & $[-1, -1, -1, -1, -1, -1, -1, -1]$  & No  \\
\hbox{} & $1.35712361i$  & \hbox{} & $\emptyset$ & \\
\hline
$(3,4,6)$ &  $x^4 + 12x^2 + 81$  & $57600$  & $[-1, -1, -1, -1]$  & No  \\
\hbox{} & $ 1.224744871 - 2.738612788i$  & \hbox{} & $\emptyset$ & \\
\hline
$(3,5,5)$ & $x^4 + 7x^2 + 11$ & $4400$ & $[-1, -1, -1, -1]$ & No \\
\hbox{} & $- 2.148961142i$ & \hbox{} & $\emptyset$ & \\
\hline
$(3,5,6)$ & $x^8 - 4x^7 + 30x^6 - 64x^5 + 262x^4 - 384x^3 + 978x^2 - 684x + 1629$  & $2734871040000$ & $[-1, -1, -1, -1, -1, -1, -1, -1]$  & No  \\
\hbox{} & $1.3660254037844386468 - 1.7150676861906891827i$ &  & $\emptyset$ &  \\
\hline
$(3,6,6)$ &  $x^2 - x + 4$ & $-15$ & $[-1, -1]$ & Yes \\
\hbox{}  & $0.500000000 + 1.936491673i$ & \hbox{} & ${\cal P}_2, {\cal P}_2'$ & \hbox{} \\
\hline
$(4,4,4)$ & $x^4 - 2x^2 - 1$ & $-1024$ & $[1, 1, -1, -1]$ & Yes \\
\hbox{} & $0.643594253i$ & \hbox{} & $\emptyset$ & \\
\hline
$(4,4,5)$ & $x^4 + 3x^2 + 1$ & $400$ & $[-1, -1, -1, -1]$ & No \\
\hbox{} & $1.618033989i$ & \hbox{} & $\emptyset$ & \\
\hline
$(4,4,6)$ & $x^4 + 2x^2 - 2$ &  $-4608$ & $[1, 1, -1, -1]$ & No* \\
\hbox{} & $1.652891650i$ & \hbox{} & $\emptyset$ &  \\
\hline
$(4,5,5)$ & $x^8 - 2x^6 - 8x^5 - 5x^4 + 8x^3 + 12x^2 + 4x - 1$ &  $368640000$ & $[1, 1, 1, 1, -1, -1, -1, -1]$ & No\\
\hbox{} & $-0.707106781 + 0.544223826i$ & \hbox{} & $$ & \\
\hline
$(4,5,6)$ & $x^8 + 2x^6 - 39x^4 - 130x^2 - 95$  & $-5042995200000$ & $[1, 1, -1, -1, -1, -1, -1, -1]$  & No \\
\hbox{} & $- 1.05321208i$ & \hbox{} & $\emptyset$  & \\
\hline
$(4,6,6)$ & $x^4 + 2x^2 + 4$ & $576$ & $[-1, -1, -1, -1]$ & No\\
\hbox{} & $0.707106781 - 1.224744871i$ & \hbox{} & $\emptyset$ &  \\
\hline
$(5,5,5)$ & $x^4 - 5$ & $-2000$ & $[1, 1, -1, -1]$ & Yes \\
\hbox{} & $1.495348781i$ & \hbox{} & $\emptyset$ & \\
\hline
$(5,5,6)$ & $x^8 + 6x^6 - 13x^4 - 66x^2 + 61$ & $12648960000$ & $[1, 1, 1, 1, -1, -1, -1, -1]$ & No \\
\hbox{} & $2.466198614i$ & \hbox{} & $\emptyset$ & \\
\hline
$(5,6,6)$ & $x^4 + 21x^2 + 99$ & $39600$ & $[-1, -1, -1, -1]$ & No\\
\hbox{} & $2.673181257i$ & \hbox{} & $\emptyset$ & \\
\hline
$(6,6,6)$ & $x^4 - 2x^3 - 2x + 1$ & $-1728$ & $[1, 1, -1, -1]$ & Yes \\
\hbox{} & $-0.3660254039 + 0.930604859i$ & \hbox{} & $\emptyset$ & \\
\hline
\end{tabular}
}
\end{center}
\caption{\label{TAB:TRUNCATED_CUBES} Commensurability invariants for some truncated cubes.}
\end{table}

\subsection{L\"obell Polyhedra}
\begin{figure}
\begin{center}
\includegraphics[scale=0.5]{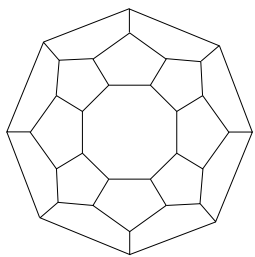}
\end{center}
\caption[]{The L\"obell polyhedron for $n=8$. \label{FIG:LOB}}
\end{figure}

For each $n \geq 5$, there is a radially-symmetric combinatorial polyhedron having two $n$-sided faces
and having $2n$ faces with $5$ sides, which provides a natural generalization of
the dodecahedron.  This combinatorial polyhedron is depicted in Figure \ref{FIG:LOB} for $n=8$.

Andreev's Theorem provides the existence of a compact right-angled polyhedron
$L_n$ realizing this abstract polyhedron because it contains no prismatic
$3$-circuits or prismatic $4$-circuits.

\begin{figure}
\begin{center}
\begin{picture}(0,0)%
\epsfig{file=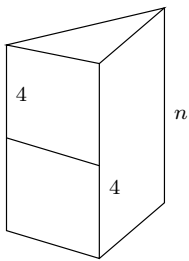}%
\end{picture}%
\setlength{\unitlength}{3947sp}%
\begingroup\makeatletter\ifx\SetFigFont\undefined%
\gdef\SetFigFont#1#2#3#4#5{%
  \reset@font\fontsize{#1}{#2pt}%
  \fontfamily{#3}\fontseries{#4}\fontshape{#5}%
  \selectfont}%
\fi\endgroup%
\begin{picture}(1246,1599)(1564,-1873)
\put(2626,-986){\makebox(0,0)[lb]{\smash{{\SetFigFont{8}{9.6}{\familydefault}{\mddefault}{\updefault}{\color[rgb]{0,0,0}$n$}%
}}}}
\put(1634,-869){\makebox(0,0)[lb]{\smash{{\SetFigFont{8}{9.6}{\familydefault}{\mddefault}{\updefault}{\color[rgb]{0,0,0}$4$}%
}}}}
\put(2218,-1453){\makebox(0,0)[lb]{\smash{{\SetFigFont{8}{9.6}{\familydefault}{\mddefault}{\updefault}{\color[rgb]{0,0,0}$4$}%
}}}}
\end{picture}%

\end{center}
\caption{The ``hexahedron'' $H_n$.  Edges labeled by an integer $n$ are
assigned dihedral angle $\frac{\pi}{n}$ and unlabeled edges are assigned
$\frac{\pi}{2}$.  \label{FIG:HEX}}
\end{figure}

An alternative construction of $L_n$ is obtained by grouping $2n$ copies of the
``hexahedron'' shown in Figure \ref{FIG:HEX} around the edge labeled $n$
\cite{VESNIN_LOB}.  This construction is shown for $L_{10}$ in Figure \ref{FIG:LOBELL10_TILING}.   We denote
this polyhedron by $H_n$, and note that, by construction, $L_n$ and $H_n$ are
commensurable for each $n$.

Of historical interest is that the first example of a closed hyperbolic
manifold was constructed by L\"obell \cite{LOBELL} in 1931 by an appropriate
gluing of $8$ copies of $L_6$.  (See also \cite{VESNIN_LOB} for an exposition
in English, and generalizations.)  This gluing corresponds to constructing an
index $8$ subgroup of the reflection group in the faces of $L_6$, hence the
L\"obell manifold has the same commensurability invariants as those presented
for $L_6$ in the table above.  It is also true that the L\"obell manifold is
arithmetic, because this underlying reflection group is arithmetic.
Arithmeticity of the classical L\"obell manifold  was previously observed by
Andrei Vesnin, but remained unpublished \cite{VESNIN_PERSONAL}.

Furthermore, Vesnin observed in \cite{VESNIN_ARITH} that the if the L\"obell
polyhedron $L_n$ is arithmetic, then $n = 5, 6, 7, 8, 10, 12,$ or $18$.  He
shows that the reflection group generated by $L_n$ contains a $(2,4,n)$
triangle group which must be arithmetic if $L_n$ is arithmetic, and applies the
classification of arithmetic triangle groups by Takeuchi \cite{TAK}.  
In combination with our computations, Vesnin's observation yields:

\begin{thm}
The L\"obell polyhedron $L_n$ is arithmetic if and only if $n = 5, 6,$ or $8$.
\end{thm}

Table \ref{TAB:LOB} contains data for the first few L\"obell polyhedra $L_n$ and, consequently, the first few ``hexahedra'' $H_n$.

\begin{table}
\begin{center}
\scalebox{0.75}{
\begin{tabular}{|l|l|l|l|l|l|}
\hline
$L_n$ & $k\Gamma$ & disc & $A\Gamma$ ramification &  Arith? \\
\hline
\hline
$L_5$ & $x^4 - x^2 - 1$ &  $-400$ & $[1, 1, -1, -1]$ & Yes \\
 & $0.786151i$ & \hbox{} & $\emptyset$ & \\
\hline
$L_6$ &  $x^2 + 2$ & $-8$ & $[-1, -1]$ & Yes \\
 & $1.414214$ & \hbox{} & ${\cal P}_3, {\cal P}_3'$ &  \\
\hline
$L_7$ & $x^6 + x^4 - 2x^2 - 1$ & $153664$ & $[1, 1, -1, -1, -1, -1]$ & No \\
\hbox{} & $-0.66711i$ & \hbox{} & $\emptyset$ & \\
\hline
$L_8$ & $x^4 - 2x^2 - 1$ &  $-1024$ & $[1, 1, -1, -1]$ & Yes \\
\hbox{} & $0.643594i$ & \hbox{} & $\emptyset$ & \\
\hline
$L_9$ & $x^6 + 3x^4 - 3$ & $1259712$ & $[1, 1, -1, -1, -1, -1]$ & No \\
\hbox{} & $1.591254i$ & \hbox{} & $\emptyset$ & \\
\hline
$L_{10}$ & $x^4 + 3x^2 + 1$ &  $400$ & $[-1, -1, -1, -1]$ & No \\
\hbox{} & $0.618034i$  & \hbox{} & $\emptyset$ & \\
\hline
$L_{11}$ & $x^{10} + 4x^8 + 2x^6 - 5x^4 - 2x^2 + 1$ & $-219503494144$ & $[1, 1, 1, 1, -1, -1, -1, -1, -1, -1]$ & No \\
\hbox{} & $1.637836i$ & \hbox{} & $\emptyset$ & \\
\hline
$L_{12}$ & $x^4 + 2x^2 - 2$ &  $-4608$ & $[1, 1, -1, -1]$ & No* \\
\hbox{} & $1.652892i$  & \hbox{} & $\emptyset$ & \\
\hline
$L_{13}$ & $x^{12} + 5x^{10} + 5x^8 - 6x^6 - 7x^4 + 2x^2 + 1$ & $564668382613504$ & not computed & No\\
\hbox{} & $1.664606i$ & \hbox{} & (long computation) & \\
\hline
$L_{14}$ & $x^6 + x^4 - 2x^2 - 1$ & $153664$ & $[1, 1, -1, -1, -1, -1]$ & No \\
\hbox{} & $1.342363i$ & \hbox{} & ${\cal P}_7, {\cal P}_7'$ &  \\
\hline
$L_{15}$ & $x^8 + 5x^6 + 5x^4 - 5x^2 - 5$ & $-1620000000$ & $[1, 1, -1, -1, -1, -1, -1, -1]$ & No \\
\hbox{} & $1.681396i$ & \hbox{} & $\emptyset$ & \\
\hline
$L_{16}$ & $x^8 + 4x^6 + 2x^4 - 4x^2 - 1$ & $-1073741824$ & $[1, 1, -1, -1, -1, -1, -1, -1]$ & No \\
\hbox{} & $1.687530i$ & \hbox{} & $\emptyset$ & \\
\hline
$L_{17}$ & not computed (long computation.) & & & \\
\hline
$L_{18}$ & $x^6 - 3x^2 - 1$ & $419904$ & $[1, 1, -1, -1, -1, -1]$ & No \\
 & $0.589319i$ & & ${\cal P}_3, {\cal P}_3'$ & \\
\hline
\end{tabular}
}
\end{center}
\caption{\label{TAB:LOB}Commensurability invariants for the L\"obell polyhedra.  Note that $L_5$ is the right-angled
	regular dodecahedron.}
\end{table}

It is interesting to notice that $L_7$ and $L_{14}$ have isomorphic invariant
trace fields, which are actually not the same field (one can check that the
specified roots generate different fields).  This alone suffices to show that
$L_7$ and $L_{14}$ are incommensurable.  Further, albeit unnecessary,
justification is provided by the fact that their invariant quaternion algebras
are not isomorphic, since $A\Gamma(L_7)$ has no finite ramification, whereas
$A\Gamma(L_{14})$ is ramified at the two finite prime ideals ${\cal P}_7, {\cal
P}_7'$ lying over the rational prime $7$.

\subsection{Modifying the L\"obell $6$ polyhedron}
Let $P$ be a compact hyperbolic polyhedron with all right dihedral angles and at
least one face $F$ with $6$ or more edges.  If one (combinatorially) splits the face $F$ into two faces $F_1$ and $F_2$ along a new edge $e$, the resulting polyhedron can be realized with all right dihedral angles so long as both of the faces $F_1$ and $F_2$ have $5$ or more edges.  The inverse of this procedure is described in \cite{INOUE_THESIS}.

\begin{figure}
\begin{center}
\begin{picture}(0,0)%
\epsfig{file=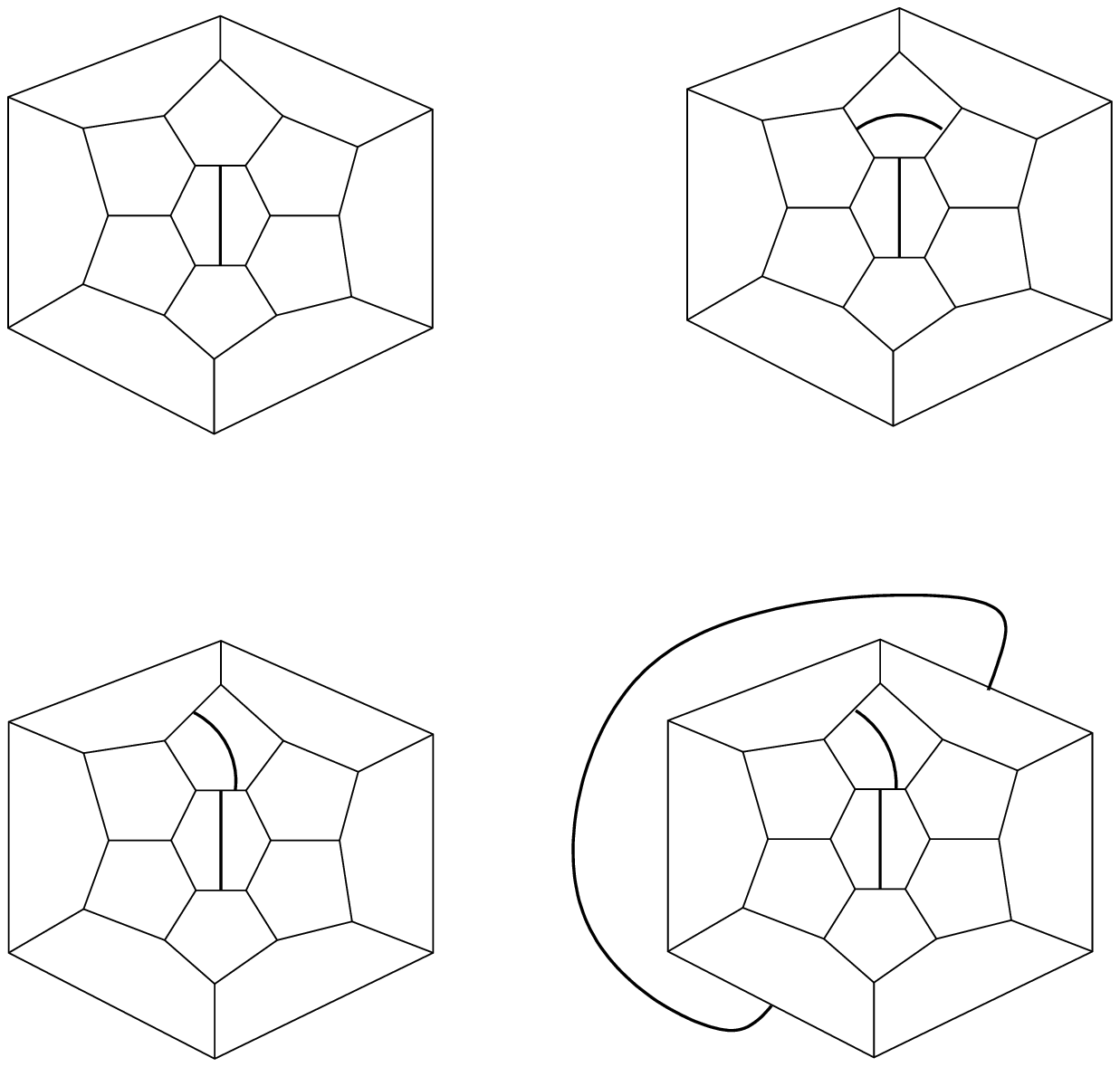}%
\end{picture}%
\setlength{\unitlength}{3947sp}%
\begingroup\makeatletter\ifx\SetFigFont\undefined%
\gdef\SetFigFont#1#2#3#4#5{%
  \reset@font\fontsize{#1}{#2pt}%
  \fontfamily{#3}\fontseries{#4}\fontshape{#5}%
  \selectfont}%
\fi\endgroup%
\begin{picture}(6064,6661)(549,-6368)
\put(3839,-5872){\makebox(0,0)[lb]{\smash{{\SetFigFont{8}{9.6}{\familydefault}{\mddefault}{\updefault}{\color[rgb]{0,0,0}$0.167273-0.768150i$}%
}}}}
\put(549,-5971){\makebox(0,0)[lb]{\smash{{\SetFigFont{8}{9.6}{\familydefault}{\mddefault}{\updefault}{\color[rgb]{0,0,0}$A\Gamma: [-1,-1,-1,-1,-1,-1,-1,-1], \emptyset$}%
}}}}
\put(549,-6121){\makebox(0,0)[lb]{\smash{{\SetFigFont{8}{9.6}{\familydefault}{\mddefault}{\updefault}{\color[rgb]{0,0,0}Not arithmetic}%
}}}}
\put(549,-5671){\makebox(0,0)[lb]{\smash{{\SetFigFont{8}{9.6}{\familydefault}{\mddefault}{\updefault}{\color[rgb]{0,0,0}$k\Gamma:  x^8 - x^6 + 3x^4 + 4x^2 + 1,  0.762867i$}%
}}}}
\put(3580,-6218){\makebox(0,0)[lb]{\smash{{\SetFigFont{8}{9.6}{\familydefault}{\mddefault}{\updefault}{\color[rgb]{0,0,0}$A_\Gamma: [1,1,-1,-1,-1,-1,-1,-1], \emptyset$}%
}}}}
\put(3580,-6368){\makebox(0,0)[lb]{\smash{{\SetFigFont{8}{9.6}{\familydefault}{\mddefault}{\updefault}{\color[rgb]{0,0,0}Not arithmetic}%
}}}}
\put(3580,-6068){\makebox(0,0)[lb]{\smash{{\SetFigFont{8}{9.6}{\familydefault}{\mddefault}{\updefault}{\color[rgb]{0,0,0}$disc = -62410000$}%
}}}}
\put(751,-2236){\makebox(0,0)[lb]{\smash{{\SetFigFont{8}{9.6}{\familydefault}{\mddefault}{\updefault}{\color[rgb]{0,0,0}$k\Gamma:  x^4 - x^3 -2x^2 - x + 1,   -0.780776 - 0.624810i$}%
}}}}
\put(751,-2386){\makebox(0,0)[lb]{\smash{{\SetFigFont{8}{9.6}{\familydefault}{\mddefault}{\updefault}{\color[rgb]{0,0,0}$disc = -1156$}%
}}}}
\put(751,-2536){\makebox(0,0)[lb]{\smash{{\SetFigFont{8}{9.6}{\familydefault}{\mddefault}{\updefault}{\color[rgb]{0,0,0}$A\Gamma:  [1,1,-1,-1],  \emptyset$}%
}}}}
\put(751,-2686){\makebox(0,0)[lb]{\smash{{\SetFigFont{8}{9.6}{\familydefault}{\mddefault}{\updefault}{\color[rgb]{0,0,0}Arithmetic}%
}}}}
\put(4276,-2386){\makebox(0,0)[lb]{\smash{{\SetFigFont{8}{9.6}{\familydefault}{\mddefault}{\updefault}{\color[rgb]{0,0,0}$disc =-7$}%
}}}}
\put(4276,-2536){\makebox(0,0)[lb]{\smash{{\SetFigFont{8}{9.6}{\familydefault}{\mddefault}{\updefault}{\color[rgb]{0,0,0}$A\Gamma: [-1,-1], {\cal P}_2, {\cal P}_2'$}%
}}}}
\put(4276,-2236){\makebox(0,0)[lb]{\smash{{\SetFigFont{8}{9.6}{\familydefault}{\mddefault}{\updefault}{\color[rgb]{0,0,0}$k\Gamma: x^2 - x + 2, 0.500000 + 1.322876$}%
}}}}
\put(4276,-2686){\makebox(0,0)[lb]{\smash{{\SetFigFont{8}{9.6}{\familydefault}{\mddefault}{\updefault}{\color[rgb]{0,0,0}Not Arithmetic (non-integral traces)}%
}}}}
\put(3589,-5667){\makebox(0,0)[lb]{\smash{{\SetFigFont{8}{9.6}{\familydefault}{\mddefault}{\updefault}{\color[rgb]{0,0,0}$k\Gamma: x^8-3x^7+5x^6-7x^5+6x^4-9x^3+5x^2-4x+1$}%
}}}}
\put(549,-5821){\makebox(0,0)[lb]{\smash{{\SetFigFont{8}{9.6}{\familydefault}{\mddefault}{\updefault}{\color[rgb]{0,0,0}$disc = 36096064$}%
}}}}
\end{picture}%
\end{center}
\caption{\label{FIG:EDGE_SURG} Right angled polyhedra obtained by adding edges to the L\"obell 6 polyhedron.}
\end{figure}

In Figure \ref{FIG:EDGE_SURG} we illustrate various right-angled polyhedra
obtained from the L\"obell 6 polyhedron by adding such edges and their commensurability invariants.

\subsection{Truncated prisms}
\begin{figure}
\begin{center}
\begin{picture}(0,0)%
\epsfig{file=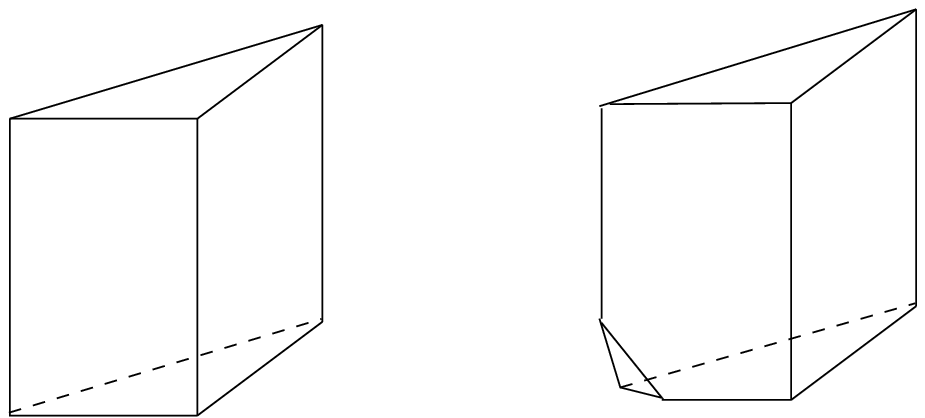}%
\end{picture}%
\setlength{\unitlength}{3947sp}%
\begingroup\makeatletter\ifx\SetFigFont\undefined%
\gdef\SetFigFont#1#2#3#4#5{%
  \reset@font\fontsize{#1}{#2pt}%
  \fontfamily{#3}\fontseries{#4}\fontshape{#5}%
  \selectfont}%
\fi\endgroup%
\begin{picture}(4587,1974)(2776,-4198)
\put(5606,-3396){\makebox(0,0)[lb]{\smash{{\SetFigFont{8}{9.6}{\familydefault}{\mddefault}{\updefault}{\color[rgb]{0,0,0}$q$}%
}}}}
\put(3941,-3461){\makebox(0,0)[lb]{\smash{{\SetFigFont{8}{9.6}{\familydefault}{\mddefault}{\updefault}{\color[rgb]{0,0,0}$3$}%
}}}}
\put(6791,-3461){\makebox(0,0)[lb]{\smash{{\SetFigFont{8}{9.6}{\familydefault}{\mddefault}{\updefault}{\color[rgb]{0,0,0}$3$}%
}}}}
\put(2776,-3496){\makebox(0,0)[lb]{\smash{{\SetFigFont{8}{9.6}{\familydefault}{\mddefault}{\updefault}{\color[rgb]{0,0,0}$q$}%
}}}}
\put(6321,-3841){\makebox(0,0)[lb]{\smash{{\SetFigFont{8}{9.6}{\familydefault}{\mddefault}{\updefault}{\color[rgb]{0,0,0}$3$}%
}}}}
\put(7106,-3991){\makebox(0,0)[lb]{\smash{{\SetFigFont{8}{9.6}{\familydefault}{\mddefault}{\updefault}{\color[rgb]{0,0,0}$3$}%
}}}}
\put(4316,-4071){\makebox(0,0)[lb]{\smash{{\SetFigFont{8}{9.6}{\familydefault}{\mddefault}{\updefault}{\color[rgb]{0,0,0}$3$}%
}}}}
\end{picture}%

\end{center}
\caption{\label{FIG:TRUNCATED_PRISMS} On the left is the family of prisms
considered in \cite{MR_BOOK}.  On the right is the family considered in this
section.  As before, edges labeled by an integer $n$ are assigned dihedral angle
$\frac{\pi}{n}$ and unlabeled edges are assigned dihedral angle
$\frac{\pi}{2}$.}
\end{figure}

Maclachlan and Reid consider triangular prisms with dihedral angles as
labeled on the left hand side of Figure \ref{FIG:TRUNCATED_PRISMS}.  Aside from the tetrahedral reflection groups and some of the Lambert cubes, this family of prisms is one of the few cases that can be computed ``by hand.''  Natural
candidates for testing our program are the truncated versions shown on the right
hand side of Figure \ref{FIG:TRUNCATED_PRISMS}.

In Table \ref{TAB:TRUNCATED_PRISMS} we show the arithmetic invariants that we have computed for the truncated prisms for $q=7,\ldots,15$.  Each of the quaternion algebras is ramified at all real places and not at any finite places, so we omit ramification data.

\begin{table}
\begin{center}\scalebox{0.75}{\begin{tabular}{|l|l|l|l|l|l|}\hline$q$ & $k\Gamma$ & disc &  Arith? \\\hline\hline
$7$ & $x^6 - 2x^5 + x^4 - 4x^3 + 3x^2 + 3x - 1$ & $-199283$ &  Yes \\
\hbox{} & $-0.400969- 1.444370i$ & \hbox{} & \\
\hline
$8$ & $x^8 - 4x^7 + 6x^6 - 7x^4 + 2x^2 + 4x - 1$ &  $473956352$ &  No \\
& $1.423880 + 1.494838i$ & \\
\hline
$9$ & $x^{12} - 4x^9 + 27x^8 - 6x^7 - 26x^6 + 6x^5 - 39x^4 + 34x^3 - 15x^2 + 6x - 1$ & $5879193047138304$ & No \\
& $-1.646799 - 1.864938i$ & \\
\hline
$10$ & $x^8 - 4x^7 + 7x^6 - 2x^5 - 5x^4 + 2x^3 + 2x^2 + 4x - 4$ & $-380000000$ & No \\
& $1.451057 + 1.553893i$ & \\
\hline
$11$ & $x^{10} - 3x^9 + x^8 + x^7 + 9x^6 - 18x^5 + 5x^4 + 9x^3 - 6x^2 - x + 1$ & $-14362045027$ & No \\
& $-1.114354 - 1.200301i$ & \\
\hline
$12$ & $x^8 - 4x^7 + 6x^6 - 4x^5 - 9x^4 + 20x^3 - 10x^2 + 1$ & $514916352$ & No \\
& $0.500000 + 1.701841i$ & \\
\hline
$13$ & $x^{12} - 5x^{11} + 10x^{10} - 17x^9 + 32x^8 - 34x^7 + 15x^6 - 11x^5 + 2x^4 + 22x^3 - 9x^2 - 6x + 1$ &  $-214921388792591$ & No \\
& $-0.470942 + 1.596598i$ & \\
\hline
$14$ & $x^{12} - 6x^{11} + 16x^{10} - 18x^9 - 2x^8 + 22x^7 - 8x^6 - 8x^5 - 16x^4 + 18x^3 + 9x^2 - 8x + 1$ &  $-581980365811712$ & No \\
& $1.474928 - 1.605111i$ & \\
\hline
$15$ & $x^8 - 2x^7 - 2x^6 + 6x^5 - 5x^4 - 4x^3 + 3x^2 + 3x + 1$ & $183515625$ & No \\
& $ -0.309017 - 0.336995i$ & \\
\hline
\end{tabular}
}
\end{center}
\caption{\label{TAB:TRUNCATED_PRISMS} Commensurability invariants for truncated prisms.}
\end{table}

\subsection{Doubly-truncated prisms} \label{SUBSEC:DOUBLE}
\begin{figure}
\begin{center}
\begin{picture}(0,0)%
\epsfig{file=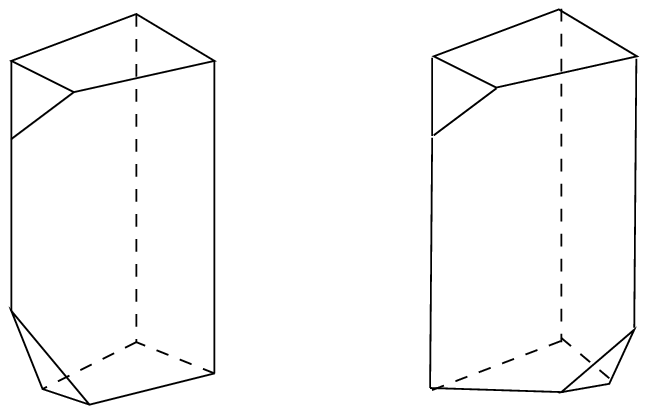}%
\end{picture}%
\setlength{\unitlength}{3947sp}%
\begingroup\makeatletter\ifx\SetFigFont\undefined%
\gdef\SetFigFont#1#2#3#4#5{%
  \reset@font\fontsize{#1}{#2pt}%
  \fontfamily{#3}\fontseries{#4}\fontshape{#5}%
  \selectfont}%
\fi\endgroup%
\begin{picture}(3394,2045)(2841,-3557)
\put(2841,-2606){\makebox(0,0)[lb]{\smash{{\SetFigFont{8}{9.6}{\familydefault}{\mddefault}{\updefault}{\color[rgb]{0,0,0}$q$}%
}}}}
\put(4051,-2686){\makebox(0,0)[lb]{\smash{{\SetFigFont{8}{9.6}{\familydefault}{\mddefault}{\updefault}{\color[rgb]{0,0,0}$q$}%
}}}}
\put(3456,-3316){\makebox(0,0)[lb]{\smash{{\SetFigFont{8}{9.6}{\familydefault}{\mddefault}{\updefault}{\color[rgb]{0,0,0}$3$}%
}}}}
\put(3661,-3496){\makebox(0,0)[lb]{\smash{{\SetFigFont{8}{9.6}{\familydefault}{\mddefault}{\updefault}{\color[rgb]{0,0,0}$3$}%
}}}}
\put(3641,-2641){\makebox(0,0)[lb]{\smash{{\SetFigFont{8}{9.6}{\familydefault}{\mddefault}{\updefault}{\color[rgb]{0,0,0}$q$}%
}}}}
\put(3406,-1871){\makebox(0,0)[lb]{\smash{{\SetFigFont{8}{9.6}{\familydefault}{\mddefault}{\updefault}{\color[rgb]{0,0,0}$3$}%
}}}}
\put(3146,-1661){\makebox(0,0)[lb]{\smash{{\SetFigFont{8}{9.6}{\familydefault}{\mddefault}{\updefault}{\color[rgb]{0,0,0}$3$}%
}}}}
\put(5671,-2666){\makebox(0,0)[lb]{\smash{{\SetFigFont{8}{9.6}{\familydefault}{\mddefault}{\updefault}{\color[rgb]{0,0,0}$q$}%
}}}}
\put(6051,-2651){\makebox(0,0)[lb]{\smash{{\SetFigFont{8}{9.6}{\familydefault}{\mddefault}{\updefault}{\color[rgb]{0,0,0}$q$}%
}}}}
\put(5206,-1596){\makebox(0,0)[lb]{\smash{{\SetFigFont{8}{9.6}{\familydefault}{\mddefault}{\updefault}{\color[rgb]{0,0,0}$3$}%
}}}}
\put(5436,-1846){\makebox(0,0)[lb]{\smash{{\SetFigFont{8}{9.6}{\familydefault}{\mddefault}{\updefault}{\color[rgb]{0,0,0}$3$}%
}}}}
\put(5316,-3526){\makebox(0,0)[lb]{\smash{{\SetFigFont{8}{9.6}{\familydefault}{\mddefault}{\updefault}{\color[rgb]{0,0,0}$3$}%
}}}}
\put(5726,-3161){\makebox(0,0)[lb]{\smash{{\SetFigFont{8}{9.6}{\familydefault}{\mddefault}{\updefault}{\color[rgb]{0,0,0}$3$}%
}}}}
\put(4851,-2681){\makebox(0,0)[lb]{\smash{{\SetFigFont{8}{9.6}{\familydefault}{\mddefault}{\updefault}{\color[rgb]{0,0,0}$q$}%
}}}}
\end{picture}%

\end{center}
\caption{\label{FIG:INCOMMENS_PRISMS} 
Edges labeled by an integer $n$ are assigned
dihedral angle $\frac{\pi}{n}$ and unlabeled edges are assigned dihedral angle
$\frac{\pi}{2}$.  (In particular, the triangular faces are at right angles to
each adjacent face.)}
\end{figure}

For $q = 4$ and $5$, there exist compact polyhedra realizing two doubly-truncated prisms pictured in Figure \ref{FIG:INCOMMENS_PRISMS}.

When $q=4$, the invariant trace fields are equal, both equal to
$\mathbb{Q}(a)$, with $a$ is an imaginary fourth root of $2$, and $A\Gamma_1 =
A\Gamma_2 \cong \left(\frac{-1 \mbox{  ,  }-1}{\mathbb{Q}\left(a\right)}
\right)$.  Indeed, one can show this using the fact that $\Gamma_1$ and $\Gamma_2$ contain subgroups isomorphic to $A_4$, see the comment at the end of Section \ref{SEC:ITF}.  However, $\Gamma_1$ and $\Gamma_2$ are incommensurable since
$\Gamma_1$ is arithmetic, while $\Gamma_2$ has non-integral traces.

If we repeat the calculation with $q=5$ both groups have invariant trace field
generated by an imaginary fourth root of $20$ and have isomorphic invariant
quaternion algebras.  However, both groups have non-integral traces, so neither
is arithmetic.  Thus, we cannot determine whether these two groups are
commensurable, or not.  (See Subsection \ref{SUBSEC:QUESTIONABLE_PAIRS}, below,
for similar examples of non-arithmetic polyhedra with matching pairs $(k\Gamma,A\Gamma)$.)

There is a good reason why these pairs have the same invariant trace field: The
invariant trace field for an amalgamated product of two Kleinian groups is the
compositum of the corresponding invariant trace fields.  See Theorem 5.6.1 from
\cite{MR_BOOK}.  The group on the left can be expressed as an amalgamated
product obtained by gluing the ``top half,'' a singly truncated prism to the
``bottom half,'' another singly truncated prism (that is congruent to the top
half) along a $(q,q,q)$-triangle group.  The group on the right can be
expressed in the same way, just with a different gluing along the $(q,q,q)$
triangle group.  This construction is similar to the construction of ``mutant
knots,'' for which commensurability questions are also delicate.  See page 190
of \cite{MR_BOOK}.

\subsection{Unexpected commensurable pairs}\label{SUBSEC:COMMENS_PAIRS}

As noted in Section \ref{SEC:ITF}, one way to find unexpected pairs of
commensurable groups is to verify that two groups are arithmetic, have the same
invariant trace field, and have isomorphic invariant quaternion algebras.  In
this section, we describe the commensurability classes of arithmetic reflection
groups in which we have found more than one group.  It is very interesting to
notice that Ian Agol has proven that there are a finite number of
commensurability classes of arithmetic reflection groups in dimension $3$,
\cite{AGOL}, however there is no explicit bound.

The reader may want to refer to the Arithmetic Zoo section of \cite{MR_BOOK} to
find other Kleinian groups within the same commensurability classes.

\begin{table}
\begin{center}\begin{tabular}{|l|l|l|l|l|l|}\hline $k\Gamma$ & disc &  $A\Gamma$ finite ramification \\\hline\hline
 $x^2 + 2$, root: $1.414213562i$ & $-8$ & ${\cal P}_3, {\cal P}_3'$ \\
The $(3,4,4)$ truncated cube &  & vol $\approx 1.0038410$\\
The L\"obell 6 polyhedron & & vol $\approx 6 \cdot 1.0038410$ \\
\hline
$x^4 - x^2 - 1$, root: $0.78615i$ & $-400$ & $\emptyset$ \\
Compact tetrahedron $T_2$ from p. 416 of \cite{MR_BOOK} & & vol $\approx 0.03588$ \\
Compact tetrahedron $T_4$ from p. 416 of \cite{MR_BOOK} & & vol $\approx 2\cdot 0.03588$ \\
The $(4,4,4)$ Lambert cube & & vol $\approx 15 \cdot 0.03588$ \\
The $(3,3,5)$ truncated cube & & vol $\approx 26 \cdot 0.03588$ \\
The L\"obell 5 polyhedron (dodecohedron) & & $\approx 120\cdot 0.03588$ \\
\hline
$x^4 - 2x^3 + x^2 + 2x - 1$, root $1.207106781 - 0.978318343i$ & $-448$ & $\emptyset$ \\
The prism on the left of Figure \ref{FIG:TRUNCATED_PRISMS}  with $q = 8$ & & vol $\approx 0.214425456$\\
The $(3,3,4)$ truncated cube & & vol $\approx  16 \cdot 0.214425456$\\
\hline
$x^4 - x^3 - x^2 - x + 1$, root: $-0.651387818 - 0.758744956i$ & $-507$ & $\emptyset$ \\
The $(3,3,3)$ Lambert cube & & vol $\approx 0.324423449$ \\
The $(3,6,6)$ Lambert cube & & vol $\approx \frac{5}{3} \cdot 0.324423449$ \\
\hline
$x^4 - 2x^2 - 1$, root: $0.643594253i$ & $-1024$ & $\emptyset$ \\
The $(4,4,4)$ truncated cube & & vol $\approx 1.1273816$\\
The L\"obell 8 polyhedron & & vol $\approx 8\cdot 1.1273816$ \\
\hline
$x^4 - 2x^3 - 2x + 1$, root: $ -0.366025 + 0.930605i$ & $-1728$ & $\emptyset$ \\
The $(3,4,4)$ Lambert cube & & vol $\approx 0.4506583058$\\
The $(6,6,6)$ truncated cube & &  vol $\approx 6 \cdot 0.4506583058$ \\
\hline
\end{tabular}
\end{center}
\caption{\label{TAB:COMMENS_CLASSES}Commensurablity classes of arithmetic reflection groups.}
\end{table}

It is easy to see that the dodecahedral reflection group is an index $120$ and
$60$ subgroup of the reflection groups in the tetrahedra $T_2$ and $T_4$,
respectively.  However, for most of these commensurable pairs the
commensurability is difficult to ``see'' directly by finding a bigger
polyhedron $Q$ that is tiled by reflections of each of the polyhedra within the
commensurability class.

For example, the commensurability of the right-angled dodecohedron $L_5$ with
the $(4,4,4)$-Lambert cube was first discovered using our computations and only
later did the authors find the explicit tiling of the dodecahedron by $8$
copies of the $(4,4,4)$-Lambert cube shown in Figure \ref{FIG:COMMENS}.  We
leave it as a challenge to the reader to observe the commensurability of the
$(3,3,5)$-truncated cube with the dodecahedron in a similar way.  (Note that at
least $13$ copies of the dodecahedron are required!)

\begin{figure}
\begin{center}
\includegraphics[scale=0.25]{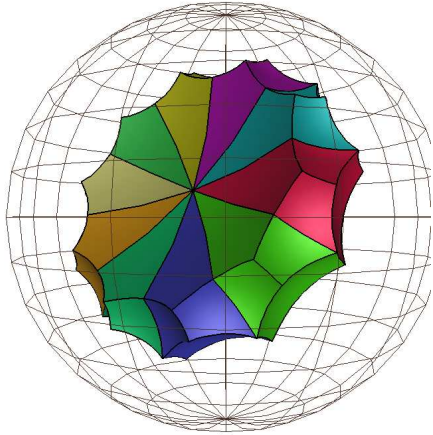}
\end{center}
\caption{\label{FIG:LOBELL10_TILING} 
Tiling the L\"obell $10$ polyhedron with $10$ copies of the $(4,4,5)$ truncated cube.}
\end{figure}

For $n=3$ and $4$ we see the commensurability between the L\"obell $2n$
polyhedron and the $(4,4,n)$ truncated cube.  This commensurability is a
general fact for any $n$: one can group $2n$ of the $(4,4,n)$ truncated cubes
around the edge with label $n$ forming the right-angled $2n$ L\"obell
polyhedron.  This construction is shown for $n=5$ in Figure \ref{FIG:LOBELL10_TILING}.

The approximate volumes were computed using Damian Heard's program Orb
\cite{ORB}.  (Orb provides twice the volume of the polyhedron, which is
the volume of the orbifold obtained by gluing the polyhedron to it's mirror
image along its boundary.)

\subsection{Pairs not distinguished by $k\Gamma$ and $A\Gamma$} \label{SUBSEC:QUESTIONABLE_PAIRS}
We found three pairs of non-arithmetic polyhedral reflection groups which are
indistinguishable by the invariant trace field and the invariant quaternion
algebra.  The first of these is the pair of
doubly truncated prisms with $q=5$, from Subsection \ref{SUBSEC:DOUBLE}.  Each of these has
the same volume, approximately $2.73952694$ since the second can be
obtained from the first by cutting in a horizontal plane and applying a
$\frac{1}{3}$ twist.
We leave it as an open question to the reader whether the two polyhedra in the
first pair are commensurable.

The second and third pairs are the $(4,4,5)$ truncated cube and the L\"obell
$10$ polyhedron and the $(4,4,6)$ truncated cube and the L\"obell $12$
polyhedron, each of which fit into the general commensurability between the
$(4,4,n)$ truncated cube and the L\"obell $2n$ polyhedron, as described in the
previous subsection.

\appendix
\section{Computing finite ramification of quaternion algebras}\label{APP:RAMCHECK}

Recall from Section \ref{SEC:QUATERNION} that to determine whether two
quaternion algebras over a number field $F$ are isomorphic it suffices to
compare their ramification over all real and finite places of $F$.

To compute the finite ramification of a quaternion algebra $A \cong
\hilb{a}{b}{F}$ we first recall that there are a finite number of candidate
primes over which $A$ can ramify: $A$ necessarily splits (is unramified) over
any prime ideal not dividing the ideal $\langle 2ab \rangle$.  To check those
primes ${\cal P}$ dividing $\langle 2ab \rangle$, we apply Theorem
\ref{THM:HILBERT_EQ} which states that $A$ splits over ${\cal P}$ if and only
if there is a solution to the Hilbert Equation  $aX^2+bY^2 = 1$ in the
completion $F_{\nu}$.  Here $\nu$ is the valuation given by $\nu(x) = c^{n_{\cal
P}(x)}$ on $R_F$, see Section \ref{SEC:QUATERNION}. 

Because it is difficult to do computer calculations in the completion
$F_{\nu}$, where elements are described by infinite sequences of elements of
$F$, we ultimately want to reduce our calculations to be entirely within
$R_F$.  Hensel's Lemma, see \cite{LANG_NUMBER}, is the standard machinery for
this reduction.  Happily, for the problem at hand, the proper use of Hensel's
Lemma has previously been worked out, see \cite{SNAP_PAPER} (whose authors write $\nu_{\cal P}$
instead of our $n_{\cal P}$ and $|\cdot|$ instead of our $\nu$).   We use the
same techniques they do but without many of the optimizations, which we find
unnecessary for our program (this is not the bottleneck in our
code).  The following theorem provides the necessary reduction:

\begin{thm}\label{THM:HENSEL_HILBERT}
Let ${\cal P}$ be a prime ideal in $R_F$, let $a,b \in R_F$ be such that
$n_{\cal P}(a), n_{\cal P}(b) \in \{0,1\}$, and define an integer $m$
as follows:  if ${\cal P} | 2$, $m=2 n_{\cal
P}(2)+3$ and if ${\cal P} \not | \, 2$ then $m=1$ if $n_{\cal P}(a) = n_{\cal
P}(b) = 0$ and $m=3$ otherwise.   

Let $S$ be a finite set of representatives for the ring $R_F/{\cal P}^m.$  The Hilbert Equation
\begin{eqnarray*}
aX^2+bY^2=1
\end{eqnarray*}
\noindent
has a solution with $X$ and $Y \in F_\nu$ if and only if there exist elements $X', Y',$ and $Z' \in S$ such that
\begin{eqnarray*}
\nu(aX'^2+bY'^2 -Z'^2) \leq c ^m
\end{eqnarray*}
and $\max\{\nu(X'), \nu(Y'), \nu(Z')\} = 1$.
\end{thm}

\noindent
Recall that $0< c < 1$ is the arbitrary constant that appears in the definition
of $\nu$.  Also note that the condition that $n_{\cal P}(a), n_{\cal P}(b) \in \{0,1\}$ is no real restriction since one can always divide the elements of $F$ appearing in the Hilbert symbol by any squares in the field without changing $A$.

Our Theorem \ref{THM:HENSEL_HILBERT} is a minor extension of Proposition 4.9 of
\cite{SNAP_PAPER}.  Their proposition only applies to dyadic primes,
corresponding to $m=2n_{\cal P}(2)+3.$   The other two cases can be proved analogously using, instead of their Lemma 4.8, the following statement:

\begin{lem}
Suppose that $\nu$ is the valuation corresponding to a non-dyadic prime ${\cal
P}$.  Let $X, X'$ be in $F_\nu$ and suppose that $\nu(X) \leq 1$ and $\nu(X-X')
\leq c^k$ for some non-negative integer $k$.  Then, $\nu(X^2-X'^2) \leq
c^{2k}$.
\end{lem}

The reason we make this minor extension of Proposition 4.9 from
\cite{SNAP_PAPER} is that we do not implement the optimizations for non-dyadic
primes that appear in SNAP, instead choosing to use the Hilbert
Equation in all cases.

%
%

Thus, the problem of determining finite ramification of $A$ reduces to finding
solutions in $R_F$ to the Hilbert inequality from Theorem
\ref{THM:HENSEL_HILBERT}.  Our program closely follows SNAP \cite{SNAP},
solving the equation by exploring, in depth-first order \cite{COMPSCI}, a tree
whose vertices are quadruples $(X,Y,Z,n)$, where $\nu(aX^2+bY^2 -Z^2) \leq c
^n$.  The children of the vertex $(X,Y,Z,n)$ are the vertices
$(X_1,Y_1,Z_1,n+1)$ where $X_1 \equiv X, Y_1 \equiv Y, \mbox{ and } Z_1 \equiv
Z \pmod{{\cal P}^n}$.   In fact, the condition $\max\{\nu(X),\nu(Y),\nu(Z)\} =
1$ means that the search space is considerably reduced because there must be a
solution to the inequality with one of $X,Y,Z$ equal to $1$.  Indeed, if
$\nu(Y) = 1$ then $Y$ is invertible modulo any power of ${\cal P}$ and
$(XY^{-1},1,ZY^{-1})$ is also a solution. So, our program only searches the
three sub-trees given by fixing $X,Y,$ or $Z$ at $1$.  In the case that we
search all three trees up to level $n=m$ unsuccessfully, then Theorem
\ref{THM:HENSEL_HILBERT} guarantees that there is no solution to the Hilbert
Equation, and consequently $A$ ramifies over ${\cal P}$.

\bibliographystyle{plain}
\bibliography{andreev.bib}
\end{document}